\documentclass[nonblindrev]{informs3}
\OneAndAHalfSpacedXI 
\usepackage{endnotes}
\let\footnote=\endnote

\usepackage[T1]{fontenc}

\usepackage{graphicx}
\usepackage{multirow}
\usepackage{hhline}
\usepackage[small, margin=1cm]{caption}
\usepackage{color}
\definecolor{strcolor}{rgb}{0.6, 0.2, 0.6}
\definecolor{commentcolor}{rgb}{0.3125, 0.5, 0.3125}
\definecolor{keycol}{rgb}{0, 0, 1}

\usepackage{amssymb}
\usepackage{amsmath}
\usepackage{mathrsfs} 
\usepackage{subfigure} 
\usepackage[linesnumbered, ruled, vlined]{algorithm2e} 
\usepackage{makecell}
\usepackage[figuresright]{rotating}
{\theoremstyle{THkey}\newtheorem{policy}{Policy}  
\usepackage{upgreek} 
\usepackage{soul} 
\usepackage{booktabs} 
\usepackage{lmodern} 
\usepackage{libertine} 
\usepackage{anyfontsize}
\usepackage{bbm}
\usepackage{listings}
\lstset{
	emph={ROVar, ROUn, ROVarDR, ROExpr, RONormInf, RONorm1, RONorm2,ROConstraint,ROExpect, ROSq, ROConstraintSet,ROIntVar,ROBinVar, ROInfinity,ROModel,ROVarDRArray, ROVarArray, ROMinimize,ROUnArray, ROAbs, ROPos, ROSum, int},emphstyle={\color{strcolor}\bfseries},
	keywordstyle={\color{blue}\bfseries},
	commentstyle={\color{commentcolor}},
	stringstyle={\color{strcolor}\bfseries},
	language=C++,               
	basicstyle={\ttfamily\footnotesize}, 
	numbers=left,                  
	numberstyle=\footnotesize,      
	stepnumber=1,                   
	numbersep=5pt,                
	backgroundcolor=\color{white},  
	showspaces=false,              
	showstringspaces=false,        
	showtabs=false,               
	frame=single,	                
	tabsize=2,	                	
	captionpos=b,                 
	breaklines=true,               
	breakatwhitespace=false,       
	escapeinside={\%*}{*)},        
	keywords=[1]{for, break, if, else, function}
}

\usepackage{url}
\newcommand {\bea}{\begin{eqnarray}}
\newcommand {\eea}{\end{eqnarray}}

\def\blot{\quad \mbox{$\vcenter{ \vbox{ \hrule height.4pt
				\hbox{\vrule width.4pt height.9ex \kern.9ex \vrule width.4pt}
				\hrule height.4pt}}$}}

\usepackage{natbib}
\bibpunct[, ]{(}{)}{,}{a}{}{,}
\def\bibfont{\fontsize{8}{9.5}\selectfont}

\TheoremsNumberedThrough    
\ECRepeatTheorems
\EquationsNumberedThrough   
\gdef\AQ#1{}
\gdef\CQ#1{}
\begin{document}
	
\def\COPYRIGHTHOLDER{INFORMS}
\def\COPYRIGHTYEAR{2017}
\def\DOI{\fontsize{7.5}{9.5}\selectfont\sf\bfseries\noindent https://doi.org/10.1287/opre.2017.1714\CQ{Word count = 9740}}

	\RUNAUTHOR{Jia et~al.} 

	\RUNTITLE{ScenPO for Data-Driven Online Inventory Routing}

\TITLE{Scenario Predict-then-Optimize for Data-Driven Online Inventory Routing}

	\ARTICLEAUTHORS{

\AUTHOR{Menglei Jia}
\AFF{Department of Management Science, Antai College of Economics and Management, Shanghai Jiao Tong University, China}

\AUTHOR{Albert H. Schrotenboer}
\AFF{Operations, Planning, Accounting, and Control Group, School of Engineering, Eindhoven University of Technology, the Netherlands}

\AUTHOR{Feng Chen}
\AFF{Sino-US Global Logistics Institute, Antai College of Economics and Management, Shanghai Jiao Tong University, China \{fchen@sjtu.edu.cn\}}
}

	\ABSTRACT{\\ \textbf{Problem definition:} The real-time joint optimization of inventory replenishment and vehicle routing is essential for cost-efficiently operating one-warehouse, multiple-retailer systems. This is complex, as future demand predictions should capture (auto)correlation and lumpy retailer demand, and based upon such predictions, inventory replenishment and vehicle routing decisions must be taken. Traditionally, such decisions are made by either making distributional assumptions or using machine-learning-based point forecasts. The former approach ignores nonstationary demand patterns, while the latter approach only provides a point forecast ignoring the inherent forecast error. Consequently, in practice, service levels often do not meet their targets, and truck fill rates fall short, harming the efficiency and sustainability of daily operations. \textbf{Methodology/results:} We propose Scenario Predict-then-Optimize. This fully data-driven approach for online inventory routing consists of two subsequent steps at each real-time decision epoch. The \textit{scenario-predict} step exploits neural networks --- specifically multi-horizon quantile recurrent neural networks --- to predict future demand quantiles, upon which we design a scenario sampling approach. The subsequent \textit{scenario-optimize} step then solves a scenario-based two-stage stochastic programming approximation. Results show that our approach outperforms a classic prediction-focused learning approach, distributional approaches, empirical sampling methods, residuals-based sample average approximation, and a state-of-the-art decision-focused learning approach. We show this both on synthetic data and large-scale real-life data from our industry partner. \textbf{Implications:} Our approach is appealing to practitioners. It is fast, does not rely on any distributional assumption, and does not face the burden of single-scenario forecasts. It also outperforms residuals-based scenario generation techniques. We show it is robust for various demand and cost parameters, enhancing the efficiency and sustainability of daily inventory replenishment and truck routing decisions. Finally, scenario Predict-then-Optimize is general and can be easily extended to account for other operational constraints, making it a useful tool in practice.}

\maketitle

\section{Introduction} \label{1}
Optimal inventory replenishment and truck routing in one-warehouse, multiple-retailer systems is fundamental for supply chain management \citep{bertazzi2008analysis}. In 2021, total operating expenses and end-of-year inventories of U.S. retailer businesses were \$1,417,106 billion and \$655.923 billion, respectively, accounting for 21.73\% and 10.06\% of total sales \citep{ARTS}. Traditionally, inventory replenishment and truck routing are controlled separately. Replenishment follows a relatively static ordering logic including base-stock and order-up-to policies \citep{huh2008s, xu2011periodic}. Only after the replenishment quantities have been determined, transportation costs associated with the distribution via truck routes are minimized. Although simple and widely used in practice, this process exhibits two apparent drawbacks: The disintegration of inventory replenishment and truck routing leading to suboptimal joint decisions, and secondly, the inability to act dynamically upon demand shifts, trends, and disruptions due to the static ordering logic. However, this traditional perspective is subject to change. Widespread data availability and efficient IT systems increasingly motivate retailers to adopt a dynamic decision-making process that \textit{integrates} replenishment and truck routing decisions based upon the latest demand information. 

Real-time inventory replenishment and truck routing are studied in (variants of) the Stochastic Inventory Routing Problem \citep{coelho2014thirty}. Generally speaking, these problems consider a depot, at which capacitated trucks are stationed, and a set of geographically scattered retailers each facing stochastic demand. Using truck routes to replenish the retailers' inventory, the goal is to minimize the expected sum of inventory and transportation costs while meeting a service-level criterion at each retailer. To do so efficiently, decision-makers rely on the quality of the predictions of future retailer demand. Indeed, such demand predictions must subsequently be merged with optimization tools to transfer the prediction into actual real-time inventory replenishment and truck routing decisions. 

However, integrating prediction and optimization is complex. Existing research has focused on imposing distributional (robust) assumptions on retailer demand \citep[see, e.g.,][]{adelman2004price, kleywegt2004dynamic, cui2023inventory}. However, demand is nonstationary and (auto)correlated in practice, which invalidates making distributional assumptions. Alternatively, researchers develop machine-learning-based point forecasts of future retailer demand \citep[see, e.g., the prediction approach in][]{huber2020daily}. While circumventing the drawbacks of (static) distributional assumptions, solely basing a decision on a single point forecast leads to overly specific decisions by ignoring forecast errors. 

As a solution, we introduce ``Scenario Predict-then-Optimize" (ScenPO), a novel data-driven prediction and optimization approach for integrated, data-driven inventory replenishment and truck-routing decisions. Our ScenPO approach neither relies on making distributional assumptions nor on a single point forecast. It consists of a distribution-free \textit{scenario-predict} step, which generates data-driven future demand scenarios for individual retailers, and a \textit{scenario-optimize} step that leverages powerful stochastic programming approaches. For the prediction step, we propose using a Multi-Horizon Quantile Recurrent Neural Network (MQRNN), which predicts future demand quantiles for a given number of periods based on past retailer demand observations \citep{wen2017multi}. We use these predicted quantiles to sample the scenarios needed for the scenario-optimize step. The scenario-optimize step then solves a multiperiod stochastic inventory routing problem based on the data-driven scenarios from the scenario-predict step. It determines joint inventory-replenishment and truck-routing decisions for a given number of future periods, of which we only execute the decisions until the next decision epoch. In this way, we do not suffer from the inherent risk of working with single-scenario forecasts. Using stochastic programming, we automatically find a solution with the lowest expected cost that hedges the forecast error. Moreover, we can leverage a suite of efficient solution approaches for such models, including scenario-decomposition methods such as progressive hedging. It is worth noting that ScenPO is generic and applies to any setting that combines demand forecasting and combinatorial optimization.

Our ScenPO approach is an alternative to so-called decision-focused learning approaches --- that integrates learning and optimization and is also called end-to-end learning --- for solving complex online combinatorial optimization problems. Instead of solely focusing on minimizing a prediction error, decision-focused learning approaches aim to minimize the expected cost of the resulting decisions \citep[see, e.g.,][]{elmachtoub2022smart, mandi2023decision}. For several reasons, such methods cannot be applied to our setting efficiently. A first stream of methods focus on the implementation of differentiable optimization layers, relying on a useful approximation of the optimization problem differential in the backward pass. Almost all existing works consider the setting where uncertainty solely lies in the objective function; however, many cases, including ours, have uncertainty in their constraints which invalidates such approaches \citep{kotary2021end}. While it is also possible to consider uncertainty in the constraints \citep[see, e.g.,][]{paulus2021comboptnet}, the learning problem has not yet been well-defined. As a consequence, it is not so effective for complex combinatorial optimization problems with dense solutions and many constraints and variables including the data-driven online inventory routing problem that we study in this paper \citep{coelho2014thirty}. A second stream develops so-called surrogate loss functions with appealing computational properties \citep[see, e.g.,][]{sahay2016multienterprise, ye2016computationally} -- which are effective in many contexts but ignore the vast stream of literature and knowledge on real-time integrated inventory replenishment and truck routing decisions.

As an efficient alternative to decision-focused learning, our ScenPO method does not learn the optimization performance but relies on stochastic programming to automatically hedge against forecast errors while considering data-driven future demand scenarios. In ScenPO, the prediction step can be performed independently of the optimization step, which is computationally appealing. We exploit stochastic programming to obtain a solution that minimizes the expected costs among the predicted scenarios. In our approach, the prediction step does not require evaluating optimization methods. Moreover, as our approach is derivative-free on the optimization problem, we also eradicate the need to develop surrogate models with favorable properties. Our ScenPO approach can thus be classified as a prediction-focused learning approach, though we reserve this term in our paper for approaches considering point forecasts to enhance readability. To the best of the authors' knowledge, our ScenPO is the first approach that integrates demand prediction and optimization in the context of joint replenishment and truck routing in data-driven online inventory routing.

In ScenPO, the prediction-step does not follow the recent trend of developing scenarios in a residuals-based manner. For example,  \citet{kannan2022data} first train a prediction model to obtain a point prediction, and then the residuals obtained during training are scaled and added to that prediction to construct scenarios. The predictors (residuals) can be of different designs to improve performance. Our ScenPO approach, and specifically the MQRNN we use, is an alternative approach, as it is based on forecasting demand quantiles which provide a robust characterization of the distribution of stochastic variables. This in turn allows us to generate demand forecasts scenarios in a robust way. 

We propose several variations of ScenPO to test and assess its different components structurally. First, besides considering MQRNNs, we also consider Long Short-Term Memory Neural Networks (LSTM) in our scenario-predict step and provide an associated scenario-sampling strategy.  For the scenario-optimize step, we consider two distinct stochastic programming formulations. First,  a single-stage (\textit{SS}) formulation assumes that the replenishment and truck-routing decisions are made in the first stage, while the second stage only evaluates the proposed solutions over the scenarios. Second, a two-stage (\textit{TS}) formulation allows for different replenishment and truck-routing decisions in all future periods. Both formulations employ a tailored matheuristic inspired by \cite{solyali2022effective} to solve the stochastic optimization problems. For the second formulation, we also employ progressive hedging \citep{rockafellar1991scenarios} to enhance its computational tractability further.

We evaluate our ScenPO approach (and its variations) compared to various benchmarks, both on the optimization and on the prediction side. That is, we consider various scenario sampling techniques that afterward use our scenario-optimize step. These are (i) a maximum likelihood approach that fits the best distribution, (ii) empirical sampling approaches, (iii) various standard prediction-focused learning approaches utilizing point forecasts, (iv) various state-of-the-art residuals-based data-driven sample average approximation approaches proposed by \citet{kannan2022data}, (v) a traditional Holt-Winters exponential smoothing approach, and (vi) a standard quantile regression model. In addition, we compare ScenPO to a state-of-the-art decision-focused learning approach called CombOptNet, as proposed by \citet{paulus2021comboptnet}. First,  we assess all the methods on synthetic data to control the degree of demand nonstationarity. Afterward, we utilize large-scale real-life data regarding the aftersales of spare parts from our industry partner, SAIC Volkswagen Automotive Co., Ltd., in China in 2020. 

Our results show that ScenPO using MQRNN combined with the \textit{SS} stochastic programming formulation performs robustly over all our instances and significantly outperforms all other methods. Our best ScenPO approach reduces the gap between an oracle solver with perfect demand knowledge and the expected value solution by 61.60\% on average, whereas the best ``classic" prediction-focused learning paradigm only reduces this gap by 30.20\% on average. Notably, the state-of-the-art sample average approximation approach by \cite{kannan2022data} reduces this gap by 51.00\%, whereas the standard quantile regression reduces this gap by 55.60\% on average. More specifically, we observe that the residuals-based approach can match the performance of our ScenPO for simple data patterns that are relatively easy to predict. However, when the data patterns become more complex, especially as observed in our real-life data, the residuals-based approach cannot match our ScenPO (and other quantile-based approaches). This confirms that quantile regression is a robust basis for scenario generation within our ScenPO framework in the context of online data-driven inventory routing.

The main contributions are summarized as follows:
\begin{itemize}
    \item We introduce Scenario Predict-then-Optimize for joint distribution and replenishment in one-warehouse-multiple-retailer systems, also known as data-driven online inventory routing. It combines prediction and optimization, provides better operational decisions, and can be easily implemented in practice. It consists of a \textit{scenario-predict} step, which provides distribution-free data-driven scenarios, which are then exploited by the \textit{scenario-optimize} step. Here, stochastic programming is employed to make joint replenishment and distribution decisions.
    \item In the context of data-driven online inventory routing, this is the first fully integrated prediction and optimization approach that does not make distributional assumptions on the input data nor works with single-point forecasts. 
    \item We detail how to employ Multi-Horizon Quantile Recurrent Neural Networks as part of our scenario-predict step, and we provide two distinct two-stage stochastic programming formulations for the scenario-optimize step. For the latter part, we provide tailored matheuristic algorithms and utilize progressive hedging to further enhance the computational performance of ScenPO. 
    \item  Both on synthetic data and large-scale real-life data from SAIC Volkswagen Automotive Co., Ltd., we show that ScenPO outperforms a range of benchmark approaches, including single-scenario prediction-focused learning, empirical sampling techniques, state-of-the-art residuals-based data-driven sample average approximation, and a state-of-the-art decision-focused learning approach using CombOptNet. 
    \item We evaluate the strengths and weaknesses of different state-of-the-art approaches through extensive numerical experiments, providing a better understanding of which paradigms are suitable for particular use cases. Our extensive numerical experiments demonstrate the value of using quantile-based scenario generation approaches over residuals-based scenario-generation approaches within ScenPO. To the best of our knowledge, this has not yet been studied in the inventory routing literature and is important to consider in the development of online data-driven algorithms.
\end{itemize}

The remainder of the paper is organized as follows. In Section~2, we review the related literature. The problem description and associated Markov decision process formulation can be found in Section~3. Subsequently, Section~4 details our ScenPO approach, and Section~5 provides several algorithms for enhancing the efficiency of the scenario-optimize step. Computational experiments on synthetic and large-scale real-life data are provided in Section~6. We conclude our work and discuss future research avenues in Section~7.

\section{Literature Review} \label{2}
The real-time joint distribution and replenishment of inventory in one-warehouse-multiple-retailer systems is often referred to as the Stochastic and Dynamic Inventory Routing Problem. In this work, we consider a new variant of this problem, which we call the Data-driven Online Inventory Routing Problem, which does not consider any distributional assumptions on retailer demands. Furthermore, we use the phrases ``warehouse" and  ``retailer" in our paper, but this can also be interpreted as ``depot" and ``customer",  which are more commonly used in the vehicle routing literature.

The classic, deterministic inventory routing problem (IRP) integrates vehicle routing and inventory management. It is motivated by vendor-managed inventory settings and has a wide range of real-life applications, from retailer replenishment to aftersales logistics of spare parts \citep{andersson2010industrial}. The Stochastic Inventory Routing Problem (SIRP) refers to the problem with stochastic demand \citep[see, e.g.,][]{sonntag2023stochastic, hasturk2024stochastic}, and the Stochastic and Dynamic Inventory Routing Problem (SDIRP) refers to the problem in which information is revealed dynamically and decisions are made online. In our review, we focus solely on the SDIRP. For a more general, routing-oriented review, we refer the reader to the excellent review by \citet{coelho2014thirty}.

The SDIRP has been studied extensively. However, most studies make distributional assumptions for retailer demand. \citet{jaillet2002delivery} concern the SDIRP in which a commodity has to be distributed to retailers repeatedly over a long time horizon, which they solve using incremental cost approximations. \citet{kleywegt2002stochastic} formulate the inventory routing problem with direct deliveries as a Markov decision process and use dynamic programming techniques to approximate a decomposed value function. That has been extended by  \citet{kleywegt2004dynamic}, in which multiple retailers are considered in a single-vehicle route. \citet{adelman2004price} also formulate the SDIRP as a Markov decision process and solve it by approximating the future costs of current actions using optimal dual prices of a linear program. \citet{coelho2014heuristics} are the first considering SDIRP without distributional assumptions, where two different policies, proactive and reactive, are proposed to deal with dynamic demand. In the proactive policy, they make decisions based on an oracle forecast, solving a deterministic equivalent formulation. Under the reactive policy, they generate decisions based on the ($s$, $S$) policy, where $s$ is calculated according to the demand distribution and can be updated at the beginning of each period to deal with dynamic demand. \citet{achamrah2022solving} adopt a similar approach to the reactive policy where they consider transshipment and substitution to mitigate shortages. \citet{crama2018stochastic} similarly consider different demand distributions over time and develop different solution methods for SDIRP for perishable products. Except for the proactive policy of \citet{coelho2014heuristics}, all other works, despite considering dynamic demand, still rely on distributional assumptions and solve the dynamic problem by considering a newly fitted distribution on a rolling basis. The most recent work is by \citet{cui2023inventory}, which employs a distributional robust optimization approach with a moment-based ambiguity set. However, although they no longer specify the distribution, corresponding moment information still needs specification. Moreover, by exploring (distributionally) robust optimization, they minimize the expected cost for the worst-case distribution within the ambiguity set and not, as in our ScenPO, minimize the expected cost.  

Slightly different works are \citet{hvattum2009using} and \citet{hvattum2009scenario}. They approximate the infinite horizon Markov decision process by solving a finite scenario tree problem on a rolling basis. However, scenario tree generation still requires a distributional assumption — they are generated by random sampling from a distribution. Different distribution assumptions for different periods can formulate dynamic demand. However, determining the underlying distributions for each future period is a fundamental yet tricky issue. The mentioned works either specify the same distribution for each period or randomly generate distribution parameters, neglecting the sequential relationship among the demand in different periods as observed in practice. The work of \citet{greif2024combinatorial} is notable, but differs from our work on a few important aspects. They solve the SDIRP through an end-to-end learning idea by adopting the ideas of \cite{parmentier2021learning, dalle2022learning, parmentier2022learning}, and \cite{baty2024combinatorial}, in which dynamic routing problems are transferred into single-period price-collecting routing problem, where the ``prices" are learned by a machine learning predictor. However, \cite{greif2024combinatorial} consider a single vehicle only, and more importantly, only consider replenishing retailers towards their full capacity. In that way, the inventory replenishment quantity is not actively decided upon, which significantly reduces decision complexity. Our approach allows for replenishing any amount of inventory. Indeed, our ScenPO approach leverages a data-driven approach that first captures future information based on historical demand through a \textit{scenario-predict} step without any distributional assumption and then gives an approximate solution through a \textit{scenario-optimize} step to hedge against the prediction error by minimizing the expected finite horizon cost using stochastic programming. It can be applied to any setting that combines forecasting and combinatorial optimization and is a generic and relatively easy to apply.

Effective solution approaches for real-life stochastic problems involve prediction and optimization components. However, most studies treat prediction and optimization independently. For example, \citet{keskin2023dynamic} adopt the classic prediction-focused learning framework to solve the waste collection problem with ``touting" as a demand management technique, where retailers who have not yet placed an order can be actively encouraged to order a service sooner. The candidate touting retailers, i.e., the retailers likely to order soon, are determined based on the predicted fill rates by a forecasting model. Then, a multi-period vehicle routing problem with touting is solved based on the predicted results. In recent years, an emerging topic has been decision-focused learning of prediction models based on optimization models. This stream of research focuses on training machine learning models to optimize the quality of the resulting decisions and, therefore, do not necessarily lead to the best point forecasts from a prediction perspective. Nevertheless, such studies are mostly limited to linear programming problems with uncertain parameters in the objective function, such as \citet{wilder2019melding}, \citet{mandi2020interior}, and \citet{elmachtoub2022smart}. \citet{liu2022online} extend the work of \citet{elmachtoub2022smart} to consider online contextual decision-making with resource constraints. However, even though their problem contains uncertainty in the constraints, the solution method is still based on the fundamental technique of \citet{elmachtoub2022smart}, transferring the constrained problem into a penalized version of the objective by introducing dual variables and using primal-dual methods. This approach applies to problems with specific structures and is not generalizable. \citet{vlastelica2019differentiation}, \citet{ferber2020mipaal}, and \citet{mandi2020smart} have extended research into combinatorial optimization, but the optimization problems targeted are still limited to problems with uncertain parameters in the objective function. Although the current state-of-the-art can also differentiate mixed integer linear programming concerning constraints \citep{paulus2021comboptnet}, the architecture needs solutions to the problem for the training process, leading to solving many optimization problems in advance. This imposes a major computational burden on the efficiency of the method. Moreover, the state-of-the-art learning problem has not yet been well-defined in this setting, and thus, is not as effective for complex combinatorial optimization problems with dense solutions and many constraints and variables. We refer the interested reader to \citet{kotary2021end, qi2022integrating, mandi2023decision, sadana2024survey} for more research in this area. 

As far as we know, we are the first to generate scenarios for a stochastic optimization problem using quantiles predicted from MQRNN. The studies of \citet{kim2023stochastic} and \citet{doan2024locational} are most relative, with quantile-based scenario generation followed by stochastic optimization. However, the forecasting method is different. \citet{kim2023stochastic} predicts quantiles using a quantile LSTM (QLSTM) \citep{wang2019probabilistic}, and \citet{doan2024locational} use a TransLSTM \citep{song2022dttrans}. Both networks generate multistep forecasts with a recursive strategy, that is, they train a model to predict the one-step-ahead estimate, and then iteratively feed this estimate back as the ground truth to forecast longer horizons. In the field of forecasting, \citet{chevillon2007direct} showed that the direct strategy, where a model directly predicts multistep forecasts, as we use within ScenPO with the MQRNN, is less biased, more stable, and more robust to model misspecification. In addition, the context of their study is also completely different, that is, they study an energy management problem that only considers day-ahead scheduling based on the forecasted energy demands and generations, which is a single-stage stochastic problem. However, the data-driven online IRP requires making decisions for multiple look-ahead periods, which is a multi-stage optimization problem and more complex to analyze.

The scenario-predict step of ScenPO is related to research on time-series forecasting. Traditional time-series forecasting methods such as exponential smoothing or the use of autoregressive models such as ARMA and its variants are widespread. In the last decade, deep learning-based forecasting has shown great promise compared to the more traditional approaches. For example, \citet{cai2019day} compare two classical deep neural network models with the seasonal ARIMAX model in day-ahead multi-step load forecasting in commercial buildings, and the best-performing deep learning model improves the forecasting accuracy by 22.6\% compared to that of the seasonal ARIMAX. In the following, we introduce the main elements of deep learning relevant for ScenPO, and refer the reader to \citet{chatfield2000time}, \citet{de200625}, and \citet{lim2021time} for more information on time-series forecasting.

Recurrent Neural Network (RNN) is a recursive neural network that takes sequences as input. It has advantages in learning the nonlinear characteristics of sequences and is widely applied in various time-series forecasting. The Long-Short Term Memory Network (LSTM) \citep{hochreiter1997long} and Gate Recurrent Unit Network (GRU) \citep{cho2014learning} are two popular variants of RNN. They addressed the issue of RNN that can only look back in time for limited timesteps due to the fed-back signal either vanishing or exploding \citep{staudemeyer2019understanding}. Both LSTM and GRU provide point forecasts, however, and by the nature of being a prediction, they are surrounded by a prediction error. In addition, most deep learning forecasting can only predict a one-step-ahead estimate and then iteratively feed this estimate back as the ground truth to forecast longer time horizons \citep{wen2017multi}. Such a process leads to prediction error accumulation - and is clearly not reliable enough to base decisions on. \citet{wen2017multi} exploit the expressiveness and temporal nature of Sequence-to-Sequence Neural Network (Seq2Seq), the nonparametric nature of Quantile Regression \citep{koenker1978regression}, which estimates the conditional quantiles of the distribution of the response variable, i.e., the $\tau$th conditional quantile of the distribution of $Y$ given $X$, $y^{\tau}$, satisfies $\mathbb{P}(Y\leq y^{\tau}|X)=\tau$, and the efficiency of Direct Multi-Horizon Forecasting, which directly predicts the multi-horizon estimates for each period instead of recursively predicting the one-step-ahead estimate. They propose a Multi-Horizon Quantile RNN (MQRNN) forecaster: a Seq2Seq framework that generates multi-horizon quantile forecasts, which does not make distributional assumptions, produces accurate probabilistic forecasts with sharp prediction intervals, and can directly predict multi-horizon with less bias, more stability, and more robustness. Therefore, we use MQRNN in our ScenPO approach to make probabilistic forecasts based on a large amount of historical data to obtain more robust results in the optimization phase. In our computational campaign, we also compare the usefulness of LSTM in the scenario-predict step of our ScenPO approach.

After the scenario-predict step, our ScenPO approach moves to the scenario-optimize step. For this, we employ two-stage stochastic programming. The Progressive Hedging Algorithm (PHA) is a horizontal decomposition method for large-scale problems proposed by \citet{rockafellar1991scenarios}, which has theoretical convergence properties when all decision variables are continuous. It can be effectively used as a heuristic when applied to problems with discrete decision variables, as in this study's case of data-driven online inventory routing. Whereas the PHA has applications in many supply chain management problems \citep{hvattum2009using, dong2015joint, kim2015optimal, hu2019multi, jalilvand2021effective, dong2022resilient}, \citet{hvattum2009using} is the only work that has explored the PHA in the context of the stochastic inventory routing. We employ the PHA for the \textit{TS} stochastic programming formulation within the scenario-optimize step of ScenPO.

\section{Problem Description} \label{3}
This section introduces the data-driven online inventory routing problem as a Markov decision process. We first introduce the main system components. We consider an infinite discrete time horizon $\mathcal K$, where each $k \in \mathcal K$ is called a period. We consider retailers $\mathcal{N}=\{1, 2, \dots, N\}$ who each face stochastic demand, and a single warehouse $\{0\}$ with unlimited stock. We denote the retailers and the warehouse by the set $\mathcal N^+$. In each period $k \in \mathcal K$, we denote the inventory at each retailer by $\mathcal{I}_k=\{I^1_k, \dots, I^N_k\}$, where $I^i_k \in \{-\infty,\ldots, I^{\max}\}$. Here, $I^{\max}$ is the inventory capacity at each retailer, and negative inventory implies backorders. Each retailer incurs per-unit holding costs $h$ for inventory that carries over between periods, and per-unit backorder costs $e$ for backlogged demand per period.
To replenish inventory from the warehouse to the retailers, we use a homogeneous set of vehicles $\mathcal V = \{1, \ldots, V\}$, where each vehicle $v \in \mathcal V$ has maximum capacity $Q$. Vehicles make, at most, a single replenishment route per period. Each replenishment route incurs variable transportation costs depending on the geographical location of the retailers.

The sequence of events in each period $k$ is as follows. First, we observe inventory levels at the retailers, then we make a vehicle routing and replenishment decision that happens instantaneously, and afterward, the demand of period $k$ occurs at each retailer. 

The goal is to find a policy that, given past demand observations and current inventory levels at the retailers, determines how much inventory to replenish to each retailer with the available fleet of vehicles that minimizes the sum of expected holding, backorder, and transportation costs per period. In the following, we formally present the formulation of the Markov decision process.
 
In any period $k$, we let $\tilde{\mathcal{Y}}_{:k} =\{\tilde{Y}_{:k}^1, \ldots, \tilde{Y}_{:k}^N\}$ denote the last $L$ demand observations at the retailers. Here, $\tilde{Y}_{:k}^i=\{\tilde y^i_{k-1},\tilde{y}^i_{k-2}, \ldots, \tilde y^i_{k-L}\}$, where $\tilde{y}^i_k$ denotes the realization of the random demand $y^i_k$ at retailer $i$ in period $k$. We assume the random variable $y^i_k$ has nonnegative discrete support $\mathscr{Y}:= \{0, 1, \ldots, O\}$. The system state $S_k$ is then denoted as $S_k=(\mathcal{I}_k, \mathcal{\tilde{Y}}_{:k})$, $S_k \in \mathcal{S}:=\left\{\boldsymbol{\times}_{i\in \mathcal{N}}\{-\infty,\ldots, I^{\max}\}\right\}\times \{\boldsymbol{\times}_{i\in \mathcal{N}}\boldsymbol{\times}_{t\in \{k-1, k-2,\dots,k-L\}}\mathscr{Y}\}$. 

After observing the system state $S_k$, the decision maker can dispatch vehicles and replenish retailers with inventory. Every decision $x_k \in \mathcal X(S_k)$ consists of a routing plan $\mathcal{R}_k=\{R_k^1, \ldots, R_k^V\}$, where we omit dependency on $x_k,S_k$ in our notation. Each $R_k^v$ consists of a set of retailers visited by vehicle $v \in \mathcal V$, or equals $\emptyset$. We visit each retailer at most once per period. 

Associated with the routing plan, we consider the replenishment plan $\mathcal{U}_k=\{U_k^1, \dots, U_k^N\}$, denoting the amount of inventory being replenished for each retailer. We ensure that the maximum inventory capacity is respected ($I_k^i + U_k^i \leq I^{\max}$) and that the vehicle capacity is respected~($\sum_{i \in R^v_k}~U_k^i~\leq~Q$). We finally impose that $U_k^i=0 \ \forall i\in \mathcal{N} \backslash \mathcal R_k$, as a retailer can only be replenished if it is visited. Formally, the decision space can thus be formulated as:
\begin{align}
\mathcal X(S_k) := \left\{ (\mathcal{R}_k, \mathcal U_k) \left |
    \begin{aligned}
    & \mathcal{R}_k {\rm \ partitions\ a\ subset\ of\ } \mathcal{N} \\
    &\begin{aligned}
    & I_k^i + U_k^i \leq I^{\max} \!\!\!\!\!&&\forall i \in \mathcal{N} \\
    & \sum_{i \in R_k^v} U_k^i \leq Q &&\forall v\in \mathcal{V}\\
    & U_k^i=0 &&\forall i\in \mathcal{N} \backslash R_k \\
    & U_k^i\geq 0 &&\forall i\in \mathcal{N}
    \end{aligned}
    \end{aligned}
    \right.
\right\} 
\end{align}

The exogenous information arriving after decision epoch $k$ is generally denoted by $W_{k+1}$ but is independent of the action $x_k$ in our setting. It consists of the demand observations in period $k$ denoted by $\tilde y_{k}^i$ for each retailer $i$.

We transition to a new state $S_{k+1}$ via a two-step transition encapsulated by the transition function $S^M(S_k, x_k, W_{k+1})$. First, we apply the replenishment plan of action $x_k$. The direct cost of decision $x_k$ is the transportation cost $\mathbb T(S_k, x_k):= \sum_{v \in V} \alpha \cdot d(R_k^v)$, where $\alpha$ is a nonnegative scaling parameter and $d(R_k^v)$ denotes the distance of the shortest vehicle route among the retailers in $R_k^v$, which can be obtained exactly by solving a Traveling Salesman Problem. The post-decision inventory levels at each retailer equal $I_k^i + U_k^i$ for all $i \in \mathcal N$. After this, we incur the retailer demand $\tilde{y}_{k}^i$ at each retailer $i$, and incur the inventory costs
\begin{align}
\mathbb{I}(S_k, x_k):= \sum_{i \in \mathcal N}\left(h\max\{0, I_k^i + U_k^i - \tilde{y}_{k}^i\} - e \min\{I_k^i + U_k^i - \tilde{y}_{k}^i, 0\}\right), 
\end{align}
where $h$ is the per unit holding cost, and $e$ is the per unit backorder cost.
We define $C(S_k, x_k) := \mathbb{T}(S_k, x_k) + \mathbb{I}(S_k, x_k)$. Then $I^i_{k+1} = I^i_k + U_k^i - \tilde{y}_{k}^i$ and $\tilde{Y}_{:k+1}^i = (\tilde{Y}_{:k}^i \backslash \{\tilde{y}_{k-L}^i\}) \cup \{\tilde{y}_{k}^i\}$ for each retailer $i$, so that $S_{k+1} = (\mathcal I_{k+1}, \mathcal{\tilde{Y}}_{:k+1})$

The goal of the data-driven online inventory routing problem is to find a policy $\pi\in\Pi$ and a decision rule $X^{\pi}: \mathcal{S}\rightarrow \mathcal{X}(S_k)$ to minimize its total expected future discounted cost:
\begin{align}
\min_{\pi \in \Pi} \mathbb{E}\left[ \sum_{k \in \mathcal K} \gamma^kC(S_k, X^\pi(S_k)) \mid  S_0 \right], \label{m3}
\end{align}

\noindent
where $S_0$ is the initial system state, $S_{k+1} = S^M(S_k, X^\pi(S_k), W_{k+1})$ and $\gamma > 0$ is a discount factor.

\section{Scenario Predict-then-Optimize}
The difficulty of solving Equation \eqref{m3} is due to the three curses of dimensionality in states, decisions, and transitions. Optimal decisions can only be obtained if each decision's future expected cost is known, commonly known to be computationally intractable. However, if a reliable description of future retailer demand can be obtained, and can be translated into a set of future demand scenarios, stochastic programming can help to determine an approximate decision that minimizes the expected cost over a smaller, finite future horizon. 

In this section, we detail our Scenario Predict-then-Optimize (ScOP) approach that leverages stochastic programming in combination with data-driven future demand scenarios. We first provide a high-level overview of the approach, after which we detail the \textit{scenario-predict} step and the \textit{scenario-optimize} step. 

\subsection{General Overview of ScenPO}
\begin{figure}[!htb]
\centering
  \includegraphics[width=0.85\linewidth]{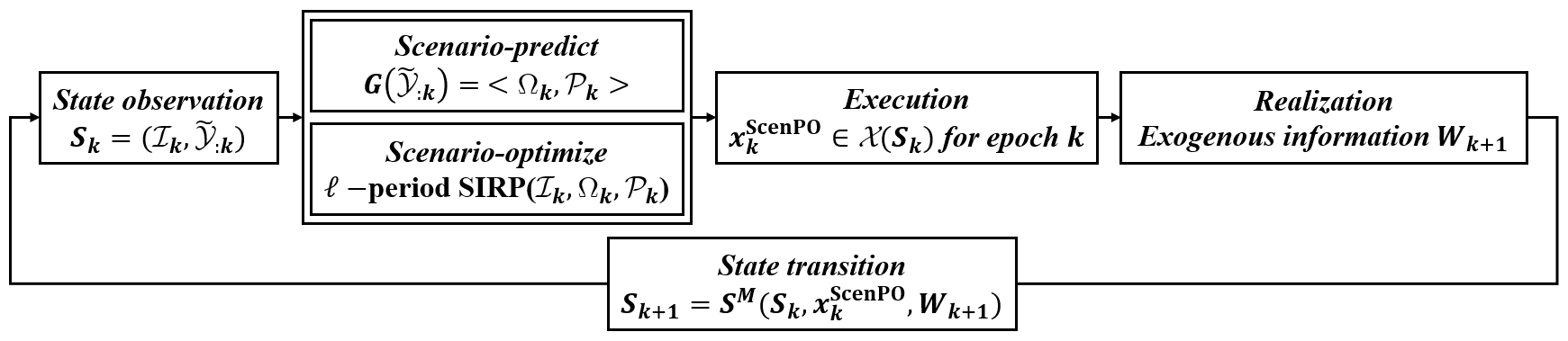}
  \caption{Graphical overview of ScenPO.}
  \label{fig_rolling_horizon}
\end{figure}

A graphical overview of ScenPO is given in Figure \ref{fig_rolling_horizon}. The inventory state transition in our data-driven online inventory routing problem equals  $I_{k+1}^i = I_k^i +U_k^i - \tilde{y}^i_{k}$. The (infinite) support of all future demand realizations equals $\left\{\boldsymbol{\times}_{i\in \mathcal{N}}\{0, \ldots, O\}\right\}^\infty$. ScenPO considers a narrowed set of future demand realizations at each decision epoch $k$ for each retailer $i$ denoted by $\Omega_k$. Here, $\omega_{k} \in \Omega_k$ describes the $\ell$ future demand realization predictions $\hat{Y}_{k:}^{i\omega_k}=\{\hat{y}_{k, k}^{i\omega_k}, \hat{y}_{k, k+1}^{i\omega_k}, \ldots, \hat{y}_{k, k+\ell-1}^{i\omega_k}\}$ for each retailer $i$, with a probability of occurring of $p(\omega_{k})$. The scenario-predict step creates a finite set $\Omega_k$ at each epoch $k$ via a machine-learning-based prediction $G$ with as input the past demand data $\mathcal{\tilde{Y}}_{:k}$, that is $G(\mathcal{\tilde{Y}}_{:k}) = \ <\Omega_k, \mathcal P_k>$, where $\mathcal P_k$ is a probability measure mapping predictions $\omega_{k}$ to $p(\omega_{k})$. The scenario-predict step thus entails finding a suitable machine learning predictor $G$ on past demand observations $\mathcal{\tilde{Y}}_{:k}$.

\begin{figure*}[!b]
    \centering
    \subfigure[Classic multistage (stochastic programming) inference]{\includegraphics[width=0.3\linewidth]{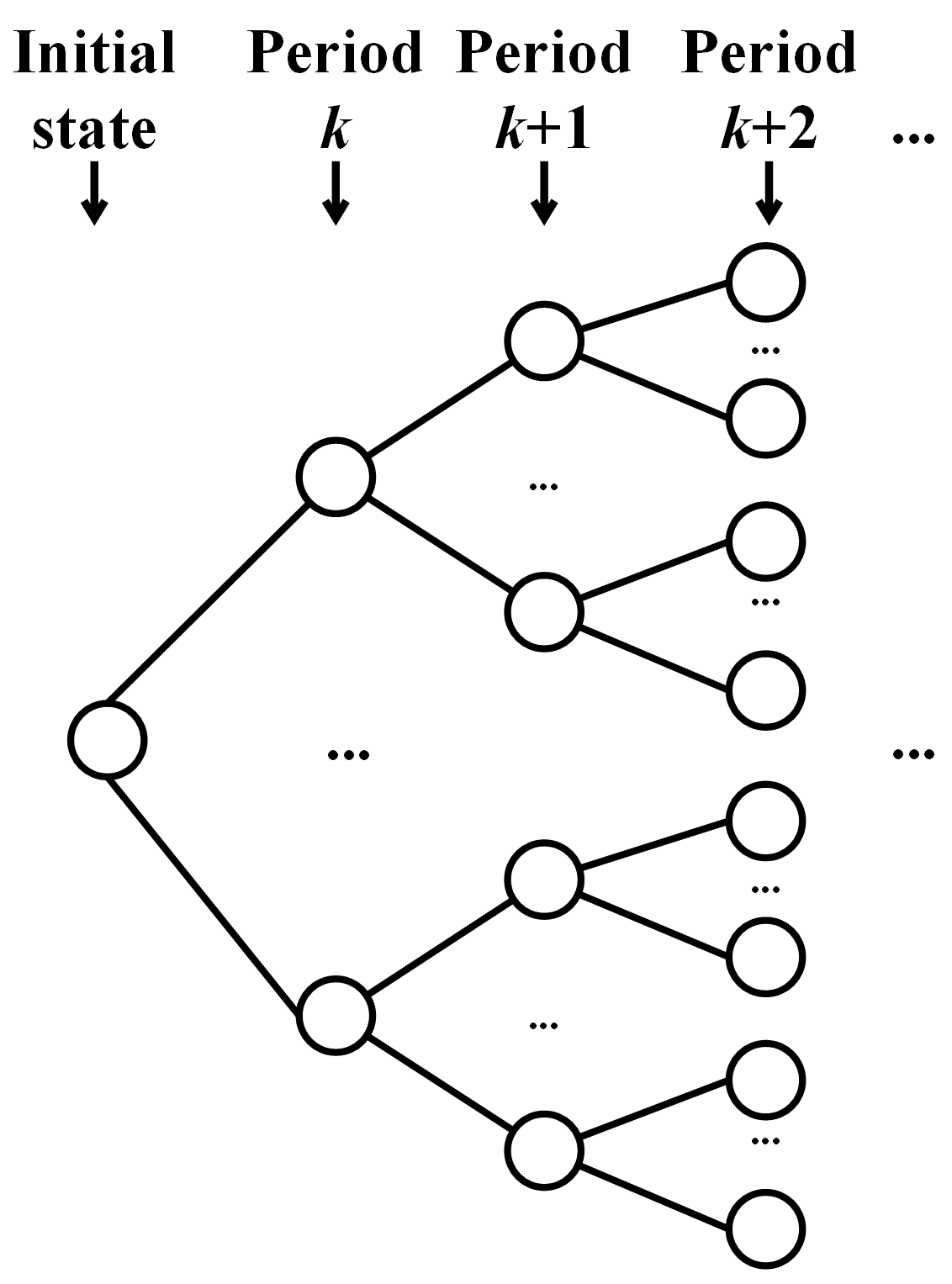}}
    \hspace{15mm}
    \subfigure[ScenPO inference]{\includegraphics[width=0.28\linewidth]{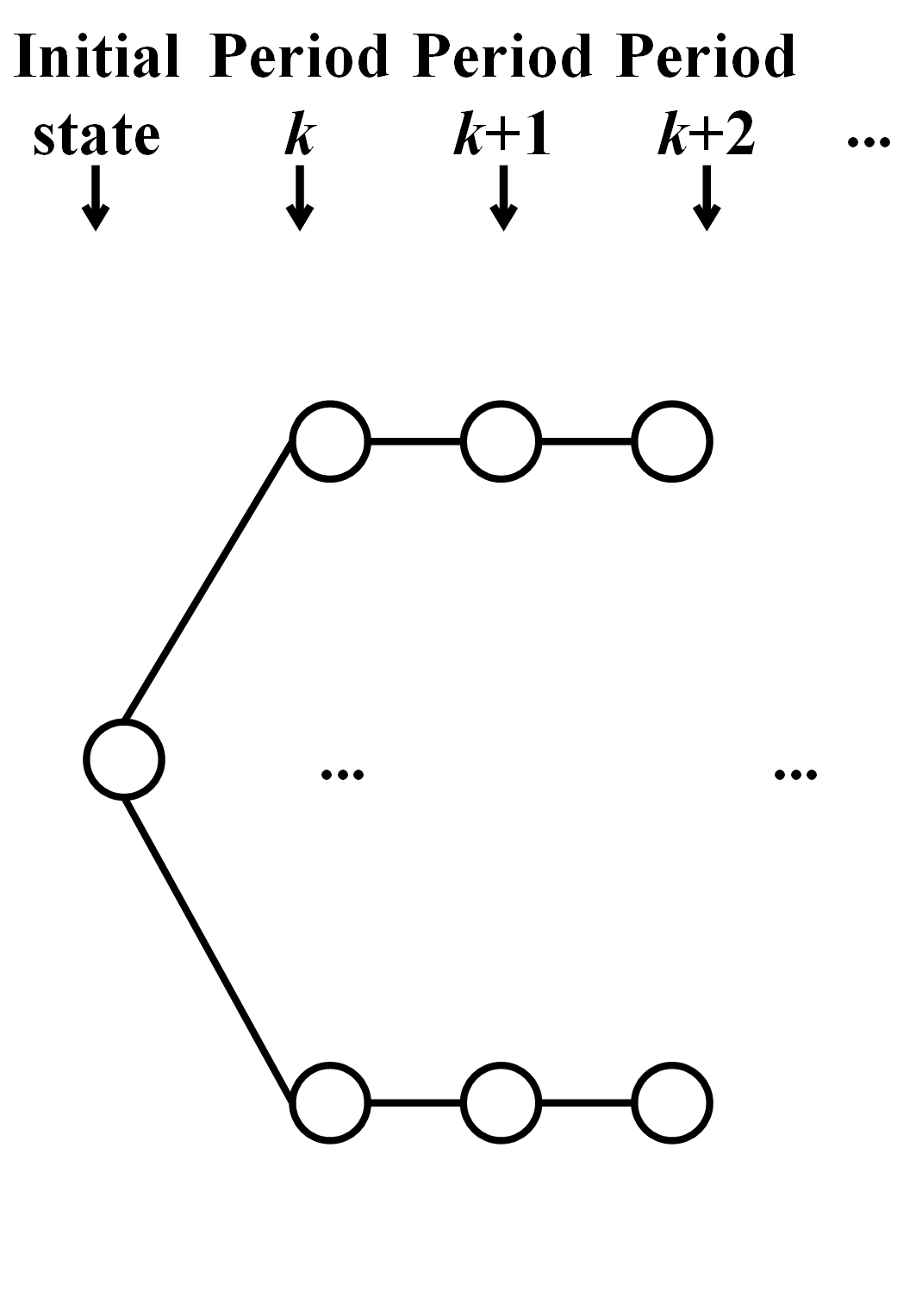}}
    \caption{The classic and ScenPO stage and decision structure.}
    \label{fig_demandRealization}
\end{figure*}

Assuming the scenario-predict step of ScenPO returns a set of future demand realizations as described by $<\Omega_k, \mathcal P_k>$, these predictions must be translated into inventory replenishment and vehicle routing decisions. The scenario-optimize step then exploits two-stage stochastic programming. This two-stage stochastic program is effectively an $\ell$-period stochastic inventory routing problem. That is, it determines for the upcoming $\ell$ periods replenishment and routing decisions that minimize the expected cost over the scenarios $\Omega_k$ obtained via predictor $G$. If $\ell$ is selected sufficiently large, the future impact of a decision in epoch $k$ can be accurately inferred if the scenarios are a valid representation of the underlying uncertainty \citep{hvattum2009scenario}. It is important to realize that within the scenario-optimize step,  our ScenPO approach considers stage aggregation \citep{powell2014clearing}. That is, we observe all future demand predictions within a scenario at once, and do not resort to multi-stage stochastic optimization. Although multi-stage stochastic optimization would give a better estimate of the impact of a decision in epoch $k$ on the future periods' decisions, it is computationally infeasible to solve $\ell$-stage stochastic optimization problems at each decision epoch. We illustrate this in Figure \ref{fig_demandRealization}. We then execute the decisions for epoch $k$ and transition to state $S_{k+1}$, after which ScenPO restarts. 

In summary, ScenPO provides a generic way to deal with MDPs for which the exogenous information exhibits nonstationarity and for which the decision space is of a combinatorial nature, as in the data-driven online inventory routing problem. We refer to the decisions of ScenPO as $x^\text{ScenPO}_k \in \mathcal X(S_k)$. In the following, we will discuss the details of the scenario-predict and the scenario-optimize step. 

\subsection{Scenario-Predict Step} \label{4}
The goal is to obtain a data-driven generator of future demand scenarios $G$. Traditional frameworks for data-driven prediction and optimization often work with point forecasts or, equivalently, mean value prediction. If the optimization is solely based upon a point forecast, it ignores the inherent forecast error of the point forecast. This only works if the point forecast is extremely accurate. However, data-driven online inventory routing and various other applications often have to deal with nonstationary demand patterns subject to exogenous disturbances, which makes it hard to provide reliable and accurate point forecasts.

The scenario-predict step of ScenPO, therefore, does not provide a single-point forecast of future retailer demands. Also, we do not rely on distributional assumptions to generate scenarios. Instead, the scenario-predict step characterizes the shape of future demand realizations by predicting future demand \textit{quantiles}. From these demand quantiles, we can then sample uniformly to effectively characterize future retailer demand scenarios without the need for distributional assumptions. 

\begin{figure*}[!htb]
    \centering
    \subfigure{\includegraphics[width=0.32\linewidth]{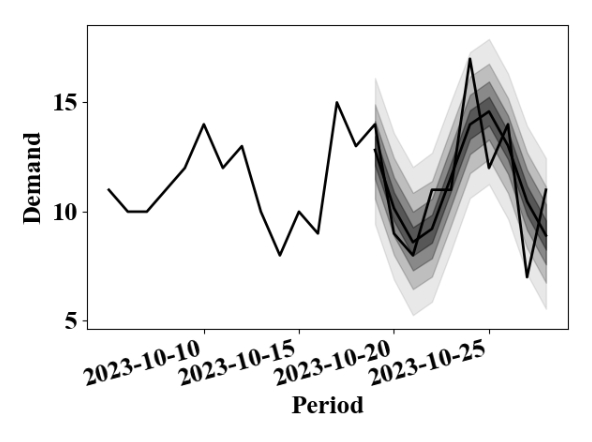}}
    \subfigure{\includegraphics[width=0.32\linewidth]{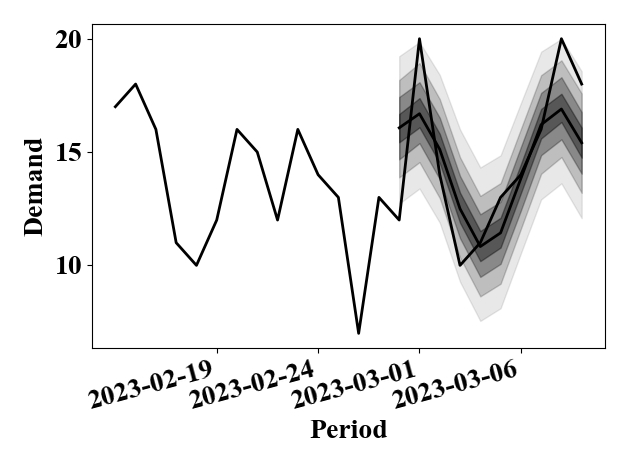}}
    \subfigure{\includegraphics[width=0.32\linewidth]{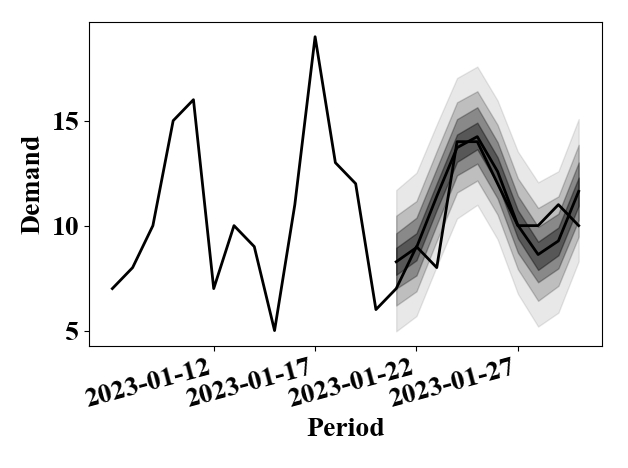}}
    \caption{Demand forecast for three example retailers by MQRNN. The different shades of gray represent the predicted quantiles.}
    \label{fig_mqrnn_forecast}
\end{figure*}

In this section, we first describe the deep learning model we adopt to predict future demand quantiles, the Multi-Horizon Quantile Recurrent Neural Network (MQRNN). Afterward, we discuss how we transform this into a scenario generation method based on the predicted quantiles by the MQRNN. An example of the performance of the MQRNN as eventually employed in our results can be found in Figure \ref{fig_mqrnn_forecast}. Here, the different shades reflect the quantile forecast of future demand.

We are interested in predicting retailer demand at each of the upcoming $\ell$ periods. We can only base our prediction on the observation of the last $L$ demand observations $\mathcal{\tilde{Y}}_{:k}=\{\tilde{Y}^1_{:k}, \dots, \tilde{Y}^N_{:k}\}$, where $\tilde{Y}^i_{:k}=\{\tilde{y}^{i}_{k-t} \mid t\in \{1, 2, \dots, L\}\}$ are the last $L$ periods demands of retailer $i$. Besides the last demand observations, we consider the periods associated with the demands as features, which can be interpreted as calendar dates or weekdays because demands often correlate with holidays or weekdays. We denote $A^l_k=\{a^{l}_{k-t} \mid t\in \{1, 2, \dots, L\}\}$ as the dates of history periods and $A^f_k=\{a^{f}_{k+t} \mid t\in \{0, 1, \dots, \ell-1\}\}$ as future dates. Using this input, the MQRNN will provide a so-called quantile forecast that gives the conditional quantiles of the target distribution, i.e., the forecasting result of a quantile $b \in \mathcal B$ and period $k$ for retailer $i$, $Y_{k:}^{\prime ib}$, means $\mathbb{P}(Y_{k:}^{i} \leq Y_{k:}^{\prime ib}|\tilde{Y}^i_{:k}, A^{l}_k, A^{f}_k)=b$. Then we have $\{ \mathcal{Y}_{k:}^{\prime b}\mid b\in \mathcal{B} \} = \mathcal{M}(\mathcal{\tilde{Y}}_{:k}, A^{l}_k, A^{f}_k)$, where $\mathcal{Y}_{k:}^{\prime b}=\{Y_{k:}^{\prime ib}\mid i\in \mathcal{N}\}$ and $\mathcal{M(\cdot)}$ is the MQRNN predictor. From this, we can then determine $<\Omega_k, \mathcal P_k>$. 

\begin{figure}[!htb]
\centering
\centering
  \includegraphics[width=0.95\linewidth]{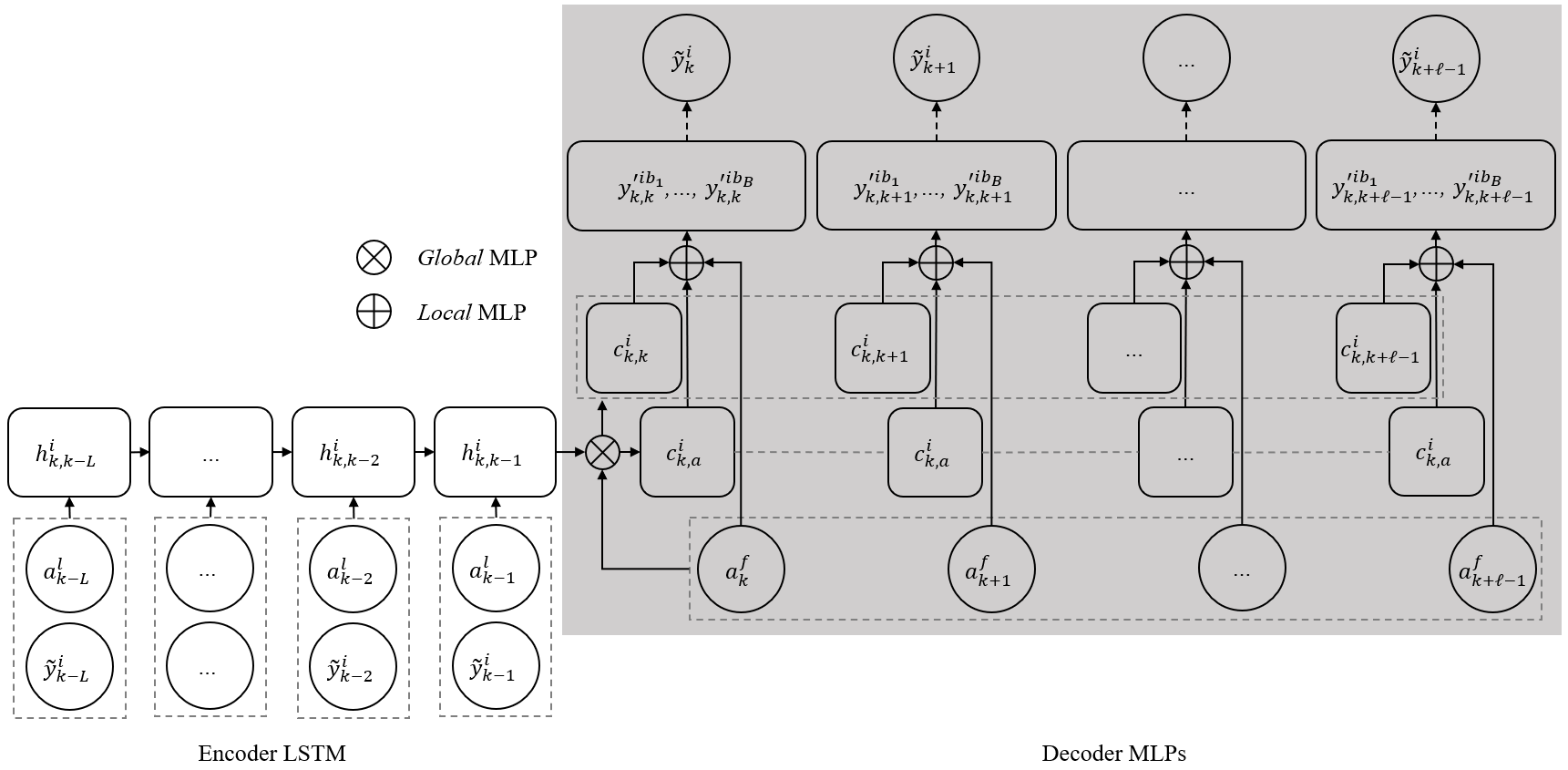}
  \caption{Visualization of the Neural Net Architecture for MQRNN for a single retailer. The dashed boxes flatten nodes into a vector, and dashed lines imply replication. The dashed arrow is the loss, which links network outputs and targets.}
  \label{fig_MQRNN_net}
\end{figure}

We now detail the MQRNN structure, which we adopt because of its demonstrated performance in demand forecasting \citep{wen2017multi, qi2023practical}. MQRNN exploits the expressiveness and temporal nature of Sequence-to-Sequence Neural Network (Seq2Seq) \citep{sutskever2014sequence}, mainly including an encoder and decoder part, see Figure \ref{fig_MQRNN_net}. MQRNN utilizes the same encoder structure as the classic Seq2Seq framework. The decoder uses a multi-horizon forecast and thus can directly train a model with a multivariate target. There are two kinds of Multi-Layer Perceptron (MLP). They are called a global MLP $m_G(\cdot)$ and a local MLP $m_L(\cdot)$. The global MLP,  $(c^i_{k,k}, c^i_{k,k+1},\dots,c^i_{k,k+\ell-1},c^i_{k, a})=m_G(h_{k,k-1}^i, A_k^f)$, summarizes the encoder output ($h^i_{k,k-1}$) plus all future inputs ($A_k^f$) into two contexts: a series of horizon-specific contexts $c^i_{k,k+t}$, where $t\in \{0, 1, \dots, \ell-1\}$, for each of the $\ell$ future horizons, and a horizon-agnostic context $c^i_{k, a}$, which captures common information. The local MLP applies to each specific horizon, which combines the corresponding future input and the two contexts from the global MLP, then outputs all the required quantiles for that specific future horizon: $(y^{\prime ib_1}_{k,k+t}, y^{\prime ib_2}_{k,k+t}, \dots, y^{\prime ib_B}_{k,k+t})=m_L(c^i_{k,k+t},c^i_{k, a}, a_{k+t}^f)$, where $\{b_1, b_2,\dots, b_B\}:=\mathcal{B}$. We feed multiple input-target pairs, i.e., ground truth demand realizations $\{\tilde y^i_{\kappa-1},\tilde{y}^i_{\kappa-2} \ldots, \tilde y^i_{\kappa-L}\}-\{\tilde y^i_{\kappa},\tilde{y}^i_{\kappa+1} \ldots, \tilde y^i_{\kappa+\ell-1}\}$ for multiple $\kappa$ and $i$, to train the encoder-decoder structure. For each input-target pair, the loss function is $\sum_{b\in \mathcal{B}} \sum_{t=0}^{\ell-1} (b(\tilde{y}^i_{\kappa+t}-y_{\kappa,\kappa+t}^{\prime ib})^+ +(1-b)(y_{\kappa,\kappa+t}^{\prime ib}-\tilde{y}^i_{\kappa+t})^+)$. The model parameters are trained to minimize the total loss of all input-target pairs. The readers are referred to \citep{wen2017multi} for more details.

Recall that we define a forecasted future demand realization of retailer $i$ in period $k$ as $\hat{Y}_{k:}^{i\omega_k}=\{\hat{y}_{k, k}^{i\omega_k}, \hat{y}_{k, k+1}^{i\omega_k}, \ldots, \hat{y}_{k, k+\ell-1}^{i\omega_k}\}$. Here we further define a scenario as $\omega_{k} = \hat{\mathcal{Y}}_{k:}^{\omega_k}=\{\hat{Y}^{i\omega_k}_{k:}\mid i\in \mathcal{N}\}$. As for the scenario generation, we first include all forecasted quantiles into our scenario set, i.e., $\Omega_k = \{ \mathcal{Y}_{k:}^{\prime b}\mid b\in \mathcal{B} \}$. Then, each time we want to generate a new scenario, we first randomly sample a quantile $b_{\beta}\in \mathcal{B}\backslash\{b_B\}$, then we construct a set of uniform distributions $\{U^i_{k,k+t}\mid i\in \mathcal{N}, t\in \{0, 1, \dots, \ell-1\}\}$, where $U^i_{k,k+t}$ is bounded by $y_{k,k+t}^{\prime ib_{\beta}}$ and $y_{k,k+t}^{\prime ib_{\beta+1}}$. Then we generate $\hat{y}_{k, k+t}^{i\omega_k}$ through randomly sampling from distribution $U^i_{k,k+t}$. After we generated all scenarios, we can easily normalize their probabilities to obtain $\mathcal P_k$.

\subsection{Scenario-Optimize Step}
The scenario-predict step searches for a decision $x^{\text{ScenPO}}_k(S_k)=x^{\text{ScenPO}}_k(\mathcal{I}_k, G(\tilde{\mathcal Y}_{:k}))$. It takes as input $G(\tilde{\mathcal Y}_{:k})$, and converts this to inventory-replenishment and vehicle-routing decisions by solving an $\ell$-period stochastic inventory routing problem.

Recall that $\mathcal{N}$ is the retailer set, and $\mathcal{N}^+$ is the set of retailers and the warehouse. Then we define the $\ell$-period stochastic inventory routing problem on an undirected graph $(\mathcal{N}^+, E)$, where $E=\{(i,j): i,j\in \mathcal{N}^+, i\neq j\}$ is the set of edges. Here, $c_{ij}$ denotes the cost of edge $(i, j)$. We let $\xi^+(i)$ and $\xi^-(i)$ be the subsets of edges starting from $i$ and ending at $i$, respectively. 

We introduce two formulations of the $\ell$-period stochastic inventory routing problem called \textit{SS} and \textit{TS}. The \textit{SS} formulation makes the $\ell$-period decisions in a single stage, i.e., makes the $\ell$-period decisions in epoch $k$ considering different scenarios of demand realizations. Formulation \textit{TS} is in a two-stage fashion: only the decisions in period $k$ are the here-and-now decisions, and the decisions for periods $k+1, k+1, \dots, k+\ell-1$ are wait-and-see decisions, which can be different according to different scenarios realization. That is illustrated in Figure \ref{fig_formulations}. 

We first introduce formulation \textit{SS}. We index the decision variables by $k$ to reflect that the system is currently in epoch $k$. Let binary variable $z^{ijv}_{k,k+t} \in \{0, 1\}$ equal to 1 if vehicle $v$ drives along edge ($i$, $j$) at period $k+t$ for $t \in \{0, \ldots, \ell - 1\}$, and 0 otherwise.  Let $u^{iv}_{k,k+t}$ denote the delivery quantity to retailer $i$ by vehicle $v$ in period $k+t$. Note that for integer demand, the delivery quantity will be integer too. Let $\mathcal{T}=\{0, 1, \dots, \ell-1\}$, and let $\hat{I}^{i\omega_k}_{k,k+t}$ denote the inventory at period $k+t$ in scenario $\omega_k \in \Omega_k$. The formulation \textit{SS} is then given by:

\begin{figure*}[!htb]
    \centering
    \subfigure[\textit{SS} illustration]{\includegraphics[width=0.43\linewidth]{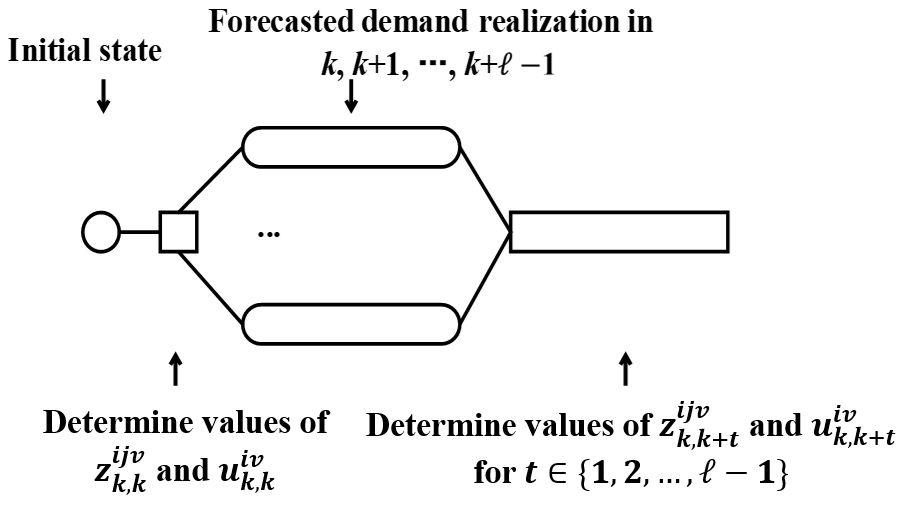}}
    \subfigure[\textit{TS} illustration]{\includegraphics[width=0.43\linewidth]{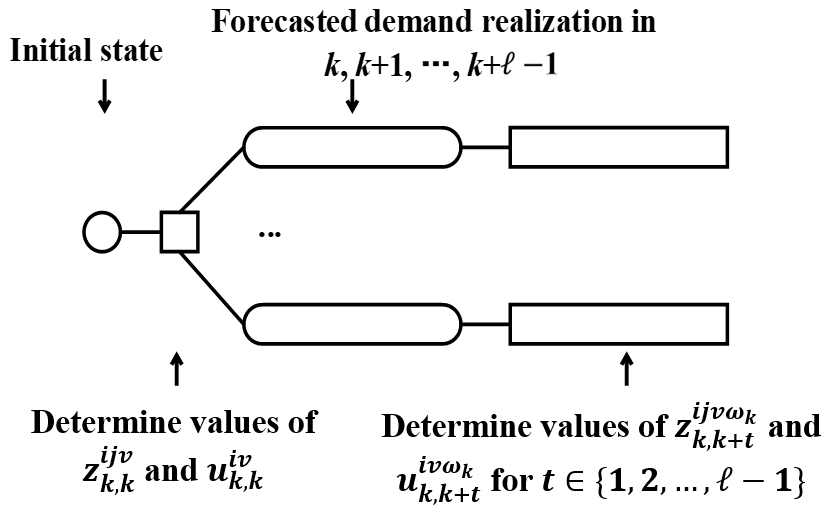}}
    \caption{Illustrations for formulations \textit{SS} and \textit{TS}.}
    \label{fig_formulations}
\end{figure*}

\begin{align}
(\textit{SS}) \ \min_{u,z}\ &  \sum_{t\in \mathcal{T}}\sum_{(i,j)\in E}\sum_{v\in \mathcal{V}}  c_{ij}z^{ijv}_{k,k+t} + \sum_{\omega_{k} \in \Omega_k} p(\omega_{k})\sum_{t\in \mathcal{T}}\!\!\!\!\!\!\!\!\!\!\!\!\!&&\sum_{i\in \mathcal{N}}(h (\hat{I}^{i\omega_k}_{k,k+t+1})^+ + e (\hat{I}^{i\omega_k}_{k,k+t+1})^-) \label{m4}\\
  s.t. &\sum_{(i,j)\in \xi^+(0)} z^{ijv}_{k,k+t}=\sum_{(i,j)\in \xi^-(0)} z^{ijv}_{k,k+t} = 1  && \forall v\in \mathcal{V}, t\in \mathcal{T}, \label{m5}\\
  &\sum_{(i,j)\in \xi^+(i)} z^{ijv}_{k,k+t} = \sum_{(j, i)\in \xi^-(i)} z^{jiv}_{k,k+t}  && \forall i\in \mathcal{N}, v \in \mathcal{V}, t\in \mathcal{T},\label{m6}\\
  &\sum_{v\in \mathcal{V}} \sum_{(i,j)\in \xi^+(i)} z^{ijv}_{k,k+t} \leq 1 && \forall i\in \mathcal{N}, t\in \mathcal{T},\label{m7}\\
  &u^{iv}_{k,k+t} \leq M\sum_{(i,j)\in \xi^+(i)} z^{ijv}_{k,k+t} && \forall i\in \mathcal{N}, v\in \mathcal{V}, t\in \mathcal{T},\label{m8}\\
  &\sum_{i\in \mathcal{N}} u^{iv}_{k,k+t} \leq Q && \forall v\in \mathcal{V}, t\in \mathcal{T},\label{m9}\\
  &u^{iv}_{k,k+t} \leq I^{\max} - \hat{I}^{i\omega_{k}}_{k,k+t} && \forall i\in \mathcal{N}, v\in \mathcal{V}, t\in \mathcal{T},\omega_{k} \in \Omega_k, \label{m10}\\
  &\hat{I}^{i\omega_{k}}_{k,k+t+1} = \hat{I}^{i\omega_{k}}_{k,k+t} + \sum_{v\in \mathcal{V}}u^{iv}_{k,k+t} - \hat{y}^{i\omega_{k}}_{k, k+t} && \forall i\in \mathcal{N}, t\in \mathcal{T}, \omega_{k} \in \Omega_k,\label{m11}\\
  &\sum_{i\in S}\sum_{j\in S, j\neq i} z^{ijv}_{k,k+t} \leq |S|-1 && \forall v\in \mathcal{V}, t\in \mathcal{T}, S\subset \mathcal{N}, 2\leq |S|\leq N,\label{m12} \\
  &z^{ijv}_{k,k+t}\in \{0,1\} && \forall (i,j)\in E, v\in \mathcal{V}, t\in \mathcal{T},\label{m13} \\
  &u^{iv}_{k,k+t}\geq 0 && \forall i\in \mathcal{N}, v\in \mathcal{V}, t\in \mathcal{T}.\label{m14}
\end{align}

The Objective Function \eqref{m4} minimizes the transportation cost and the expected $\ell$-period holding and backorder costs under different scenarios. Constraints \eqref{m5} ensure vehicles start from and end at the warehouse, while Constraints \eqref{m6} ensure flow conservation. Constraints \eqref{m7} ensure that each retailer can be visited at most once. Constraints \eqref{m8} model that a vehicle can only make deliveries to a retailer when it visits the retailer. Constraints \eqref{m9} and \eqref{m10} are vehicle capacity and maximum inventory restrictions, respectively. Constraints \eqref{m11} enforce the inventory flow balance, and Constraints \eqref{m12} are subtour elimination constraints.  Constraints \eqref{m13} and \eqref{m14} indicate the variables' domain.

Compared to formulation \textit{SS}, formulation \textit{TS} can make different routing and replenishment decisions for each scenario $\omega_k \in \Omega_k$. For this formulation, we therefore index the variables $u$ and $z$ with a scenario index $\omega_k$. 

\begin{align}
(\textit{TS}) \ \min_{u,z}\ & \sum_{\omega_{k} \in \Omega_k} p(\omega_{k})(\sum_{t\in \mathcal{T}}\sum_{(i,j)\in E}\sum_{v\in \mathcal{V}}  c_{ij}z^{ijv\omega_{k}}_{k,k+t}+\sum_{t\in \mathcal{T}}\!\!\!\!\!\!&&\sum_{i\in \mathcal{N}}(h (\hat{I}^{i\omega_{k}}_{k,k+t+1})^+ + e (\hat{I}^{i\omega_{k}}_{k,k+t+1})^-)), \label{m15}\\
  s.t. &\sum_{(i,j)\in \xi^+(0)} z^{ijv\omega_{k}}_{k,k+t}=\sum_{(i,j)\in \xi^-(0)} z^{ijv\omega_{k}}_{k,k+t} = 1 && \forall v\in \mathcal{V}, t\in \mathcal{T}, \omega_{k} \in \Omega_k,\label{m16}\\
  &\sum_{(i,j)\in \xi^+(i)} z^{ijv\omega_{k}}_{k,k+t} = \sum_{(j, i)\in \xi^-(i)} z^{jiv\omega_{k}}_{k,k+t} && \forall i\in \mathcal{N}, v \in \mathcal{V}, t\in \mathcal{T}, \omega_{k}\in\Omega_k,\label{m17}\\
  &\sum_{v\in \mathcal{V}} \sum_{(i,j)\in \xi^+(i)} z^{ijv\omega_{k}}_{k,k+t} \leq 1 && \forall i\in \mathcal{N}, t\in \mathcal{T},\omega_{k} \in \Omega_k,\label{m18}\\
  &u^{iv\omega_{k}}_{k,k+t} \leq M\sum_{(i,j)\in \xi^+(i)} z^{ijv\omega_{k}}_{k,k+t} && \forall i\in \mathcal{N}, v\in \mathcal{V}, t\in \mathcal{T},\omega_{k}\in\Omega_k, \label{m19}\\  
  &\sum_{i\in \mathcal{N}} u^{iv\omega_{k}}_{k,k+t} \leq Q && \forall v\in \mathcal{V}, t\in \mathcal{T}, \omega_{k}\in\Omega_k,\label{m20}\\
  &u^{iv\omega_{k}}_{k,k+t} \leq I^{\max} - \hat{I}^{i\omega_{k}}_{k,k+t} && \forall i\in \mathcal{N}, v\in \mathcal{V}, t\in \mathcal{T}, \omega_{k} \in \Omega_k,\label{m21}\\
  &\hat{I}^{i\omega_k}_{k,k+t+1} = \hat{I}^{i\omega_{k}}_{k,k+t} + \sum_{v\in \mathcal{V}}u^{iv\omega_{k}}_{k,k+t} - \hat{y}^{i\omega_{k}}_{k, k+t} && \forall i\in \mathcal{N}, t\in \mathcal{T}, \omega_{k} \in \Omega_k, \label{m22}\\
  &\sum_{i\in S}\sum_{j\in S, j\neq i} z^{ijv\omega_{k}}_{k,k+t} \leq |S|-1 && \forall v\in \mathcal{V}, t\in \mathcal{T}, \omega_{k}\in\Omega_k, S\subset \mathcal{N}, 2\leq|S|\leq N, \label{m23}\\
  &u^{iv\omega_{k}}_{k,k} = \bar{u}^{iv}_{k,k} && \forall i\in \mathcal{N}, v\in \mathcal{V}, \omega_{k} \in \Omega_k, \label{m24}\\
  &z^{ijv\omega_{k}}_{k,k} = \bar{z}^{ijv}_{k,k} && \forall (i,j)\in E, v\in \mathcal{V}, \omega_{k} \in \Omega_k, \label{m25}\\
  &z^{ijv\omega_{k}}_{k,k+t}\in \{0,1\} && \forall (i,j)\in E, v\in \mathcal{V}, t\in \mathcal{T},\omega_{k}\in \Omega_k, \label{m26}\\
  &u^{iv\omega_{k}}_{k,k+t}\geq 0 && \forall i\in \mathcal{N}, v\in \mathcal{V}, t\in \mathcal{T}, \omega_{k}\in\Omega_k.\label{m27}
\end{align}

The Objective Function \eqref{m15} of \textit{TS} minimizes the expected $\ell$-period transportation, holding, and backorder costs. Compared to \textit{SS}, there are two more sets of constraints \eqref{m24} and \eqref{m25}, which are nonanticipativity constraints to ensure the consistency of decisions among the scenarios in period $k$.

The optimal solution of \textit{SS} is a feasible solution of \textit{TS}. Formulation \textit{SS} is relatively conservative as we search for a single solution that performs well in all scenarios and do not employ a state-dependent evaluation of decisions in future periods. In contrast, formulation \textit{TS} is relatively optimistic as it allows for state-dependent evaluation of decisions in future periods. However, our state aggregation bases these state-dependent decisions on information that, in practice, is not yet known. Therefore, and because we employ both formulations in a rolling horizon framework, it is unclear which will perform better in practical situations. Among others, we evaluate these differences in Section \ref{6}.

After solving either formulation \textit{SS} or \textit{TS}, we let the routing plan ${R}^v_k=\{i\in \mathcal{N} \mid \sum_{j\in \mathcal{N}^+}z^{ijv}_{k,k} = 1\}$, and the replenishment plan $\mathcal{U}_k=\{\sum_{v\in\mathcal{V}}u^{iv}_{k,k}\mid i\in \mathcal{N}\}$. Then the decision for epoch $k$ can be denoted as $x^{\text{ScenPO}}_k(S_k)=(\mathcal{R}_k, \mathcal{U}_k)$.

\section{Algorithms for the Scenario-Optimize Step} \label{5}
Formulation \textit{TS}  has second-stage binary variables, making it much more complex to solve than \textit{SS}. Also, \textit{TS} requires the decision in the first period to meet nonanticipativity constraints for all upcoming periods in the aggregated stage. Such formulations clearly cannot be solved by off-the-shelf commercial MIP solvers. Therefore we propose using a Progressive Hedging Algorithm (PHA), which we will outline in Section \ref{5.1}. Following that, we present the matheuristic algorithm proposed by \citet{solyali2022effective} in Section \ref{5.2}, which we adopt to solve the \textit{SS} formulation and the subproblems arising in the PHA for the \textit{TS} formulation.

\subsection{Progressive Hedging Algorithm (PHA)} \label{5.1}
The PHA decomposes a multiple-scenario formulation into multiple single-scenario formulations by relaxing the nonanticipativity constraints and then iteratively solving penalized single-scenario formulations, known as subproblems, to reach consensus amongst the first-stage decision variables. Here, consensus means that they converge to the same value, which is forced via iteratively updating the penalization. Before we detail this, we introduce the one-scenario formulation of \textit{TS} under scenario $\omega_{k}$, called \textit{OS($\omega_{k}$)}, as follows:

\begin{align}
\textit{OS(}\omega_k\textit{)} \ \min_{u,z}\ & \sum_{t\in \mathcal{T}}\sum_{(i,j)\in E}\sum_{v\in \mathcal{V}}  c_{ij}z^{ijv\omega_{k}}_{k,k+t}+\sum_{t\in \mathcal{T}} \sum_{i\in \mathcal{N}}(h (\hat{I}^{i\omega_{k}}_{k,k+t+1})^+ \!\!\!\!\!\!\!&&+ e (\hat{I}^{i\omega_{k}}_{k,k+t+1})^-),\\
  s.t. &\sum_{(i,j)\in \xi^+(0)} z^{ijv\omega_{k}}_{k,k+t}=\sum_{(i,j)\in \xi^-(0)} z^{ijv\omega_{k}}_{k,k+t} = 1 && \forall v\in \mathcal{V}, t\in \mathcal{T},\\
  &\sum_{(i,j)\in \xi^+(i)} z^{ijv\omega_k}_{k,k+t} = \sum_{(j, i)\in \xi^-(i)} z^{jiv\omega_{k}}_{k,k+t} && \forall i\in \mathcal{N}, v \in \mathcal{V}, t\in \mathcal{T},\\
  &\sum_{v\in \mathcal{V}} \sum_{(i,j)\in \xi^+(i)} z^{ijv\omega_{k}}_{k,k+t} \leq 1 && \forall i\in \mathcal{N}, t\in \mathcal{T},\\
  &u^{iv\omega_{k}}_{k,k+t} \leq M\sum_{(i,j)\in \xi^+(i)} z^{ijv\omega_{k}}_{k,k+t} && \forall i\in \mathcal{N}, v\in \mathcal{V}, t\in \mathcal{T} ,\\  
  &\sum_{i\in \mathcal{N}} u^{iv\omega_{k}}_{k,k+t} \leq Q && \forall v\in \mathcal{V}, t\in \mathcal{T},\\
  &u^{iv\omega_{k}}_{k,k+t} \leq I^{\max} - \hat{I}^{i\omega_{k}}_{k,k+t} && \forall i\in \mathcal{N}, v\in \mathcal{V}, t\in \mathcal{T},\\
  &\hat{I}^{i\omega_k}_{k,k+t+1} = \hat{I}^{i\omega_{k}}_{k,k+t} + \sum_{v\in \mathcal{V}}u^{iv\omega_{k}}_{k,k+t} - \hat{y}^{i\omega_{k}}_{k, k+t} && \forall i\in \mathcal{N}, t\in \mathcal{T},\\
  &\sum_{i\in S}\sum_{j\in S, j\neq i} z^{ijv\omega_{k}}_{k,k+t} \leq |S|-1 && \forall v\in \mathcal{V}, t\in \mathcal{T}, S\subset \mathcal{N}, 2\leq |S|\leq N,\\
  &z^{ijv\omega_{k}}_{k,k+t}\in \{0,1\} && \forall (i,j)\in E, v\in \mathcal{V}, t\in \mathcal{T} ,\\
  &u^{iv\omega_{k}}_{k,k+t}\geq 0 && \forall i\in \mathcal{N}, v\in \mathcal{V}, t\in \mathcal{T}.
\end{align}

Clearly, we cannot just solve each subproblem independently, as it will violate the nonanticipativity constraints. To force nonanticipativity, we can add these constraints with consensus variables to each of the subproblems. Then, we use an augmented Lagrangian relaxation to relax these constraints and add a penalty term in the objective function. Thus, we consider individual scenario subproblems where the objective function can be rewritten as Equation \eqref{m39}, where $\lambda_{k}^{iv\omega_{k}}$ denotes Lagrangian multipliers and $\rho_k$ the penalty parameters. 

\begin{align}
f(\omega_k) = \min\ &\sum_{t\in \mathcal{T}}\sum_{(i,j)\in E}\sum_{v\in \mathcal{V}}  c_{ij}z^{ijv\omega_{k}}_{k,k+t}+\sum_{t\in \mathcal{T}} \sum_{i\in \mathcal{N}}(h (\hat{I}^{i\omega_{k}}_{k,k+t+1})^+ + e (\hat{I}^{i\omega_{k}}_{k,k+t+1})^-)+\nonumber\\
&\sum_{i\in \mathcal{N}}\sum_{v\in \mathcal{V}} \lambda_{k}^{iv\omega_k}(u_{k,k}^{iv\omega_{k}}-\bar{u}_{k,k}^{iv}) + \frac{1}{2}\rho_k\sum_{i\in \mathcal{N}}\sum_{v\in \mathcal{V}}(u_{k,k}^{iv\omega_{k}}-\bar{u}_{k,k}^{iv})^2,\label{m39}
\end{align}

The above objective has a quadratic penalty component containing the integer decision variable. Thus, we adopt the method by \citet{long2012sample} to force $u_{k,k}^{iv\omega_{k}}$ to converge by augmenting the objective function. The objective function can be expressed as Equation \eqref{m40}. Then, Equations \eqref{m42} and \eqref{m43} are introduced to linearize the absolute form in Equation \eqref{m40}, and the objective function is rewritten as Equation \eqref{m41}.

\begin{align}
f(\omega_k) = \min\ &\sum_{t\in \mathcal{T}}\sum_{(i,j)\in E}\sum_{v\in \mathcal{V}}  c_{ij}z^{ijv\omega_{k}}_{k,k+t}+\sum_{t\in \mathcal{T}} \sum_{i\in \mathcal{N}}(h (\hat{I}^{i\omega_{k}}_{k,k+t+1})^+ + e (\hat{I}^{i\omega_{k}}_{k,k+t+1})^-)+\nonumber\\
& \sum_{i\in \mathcal{N}}\sum_{v\in \mathcal{V}}\lambda_{k}^{iv\omega_k}|u_{k,k}^{iv\omega_{k}}-\bar{u}_{k,k}^{iv}| \label{m40}\\
f(\omega_k) = \min\ &\sum_{t\in \mathcal{T}}\sum_{(i,j)\in E}\sum_{v\in \mathcal{V}}  c_{ij}z^{ijv\omega_{k}}_{k,k+t}+\sum_{t\in \mathcal{T}} \sum_{i\in \mathcal{N}}(h (\hat{I}^{i\omega_{k}}_{k,k+t+1})^+ + e (\hat{I}^{i\omega_{k}}_{k,k+t+1})^-)+\nonumber\\
& \sum_{i\in \mathcal{N}}\sum_{v\in \mathcal{V}}\lambda_{k}^{iv\omega_k}(\mu_{k,k}^{iv} + \nu_{k,k}^{iv}) \label{m41}\\
s.t. & \ 
\mu_{k,k}^{iv} - \nu_{k,k}^{iv} = u_{k,k}^{iv\omega_{k}}-\bar{u}_{k,k}^{iv} \quad \forall i\in \mathcal{N}, v\in \mathcal{V} \label{m42}\\
& \mu_{k,k}^{iv}, \nu_{k,k}^{iv} \geq 0  \quad \quad\quad \quad \quad \quad \  \forall i\in \mathcal{N}, v\in \mathcal{V} \label{m43}
\end{align}

There are two types of variables, the routing variables $\textbf{z}$ and the replenishment variables $\textbf{u}$. Regarding the consistency of the variables, we have the following result.
\begin{proposition}[Consistency]
\label{prob:variables1}
The consistency of $\textbf{u}$ among all scenarios is sufficient for the consistency of $\textbf{z}$.
\end{proposition}

\textbf{Proof}. The only difference among scenarios is the demands faced by retailers. Thus, there is no difference among scenarios once the delivery variable of $\textbf{u}$ is consistent, which means the optimal solution of $\textbf{z}$ is consistent.\hfill \Halmos

\begin{remark}
In the case of degeneracy, we adopt the same seed to solve the subproblems, leading the subproblems to obtain the same routing decision when $\textbf{u}$ is consistent.
\end{remark}

Note that consistency of $\textbf{z}$ only implies that the maximum total delivery quantity to retailers visited by the same vehicle is consistent. Individual retailer deliveries can differ, so the optimal solution of $\textbf{u}$ can be inconsistent, given that $\textbf{z}$ is consistent.  According to Proposition \ref{prob:variables1}, once $\textbf{u}$ is consistent, the optimal solution of $\textbf{z}$ is consistent.  For that reason, we will only enforce nonanticipativity on the $\textbf{u}$ variables. Furthermore, we observe that small changes in $\textbf{u}$ will probably not lead to drastic changes to the optimal solution, whereas small changes in  $\textbf{z}$  might lead to drastically different solutions. This observation motivates us to use a myopic policy in our PHA:

\begin{policy} [Myopic Policy.] \label{policy}
During the iterations of PHA, once $\textbf{z}$ is consistent between two subsequent iterations, we fix $\textbf{z}$ in all future iterations of PHA.
\end{policy}

As for the consensus variables calculation, parameter updating, and termination criteria, we adopt a similar approach as \citet{hu2019multi}. The consensus variables are calculated as Equation \eqref{m44}. Equations \eqref{m45}-\eqref{m49} show the update of PHA's parameters, where we use superscript of ($r$), $r\geq 0$, to represent the parameters under different iterations of PHA and (0) representing the start of the algorithm. After determining the optimal routing decision of $\textbf{z}$, we set the parameters update more aggressively, as in Equation \eqref{m47}, where $\beta^0=1$ if $\textbf{z}$ is inconsistent, and $\beta^0>1$ otherwise. The termination criteria are given by Equation \eqref{m50}. 

\begin{align}
& \bar{u}_{k,k}^{iv} = \sum_{\omega_{k}\in\Omega_k} p(\omega_{k}) u_{k,k}^{iv\omega_k} \quad \forall i\in \mathcal{N}, v\in \mathcal{V}, \label{m44}\\
& \theta^{P(r)}_k = \sum_{i\in \mathcal{N}}\sum_{v\in \mathcal{V}}(\bar{u}_{k,k}^{iv(r)}-\bar{u}_{k,k}^{iv(r-1)})^2 ,\label{m45}\\
& \theta^{D(r)}_k = \sum_{i\in \mathcal{N}}\sum_{v\in \mathcal{V}}\sum_{\omega_k \in \Omega_k}(u_{k,k}^{iv\omega_{k}(r)}-\bar{u}_{k,k}^{iv(r)})^2, \label{m46}\\
& \rho^{(r)}_k=
\begin{cases}
 \beta^0 \beta^D \rho^{(r-1)}_k, &\text{if } \  \theta^{D(r-1)}_k - \theta^{D(r-2)}_k > 0\\
 \frac{1}{\beta^0\beta^P} \rho^{(r-1)}_k, &\text{if } \  \theta^{P(r-1)}_k - \theta^{P(r-2)}_k > 0\\
 \rho^{(r-1)}_k,&\text{otherwise}
\end{cases}
 \label{m47}\\
& \lambda_{k}^{iv\omega_k(0)} = \rho^{(0)}_k|u_{k,k}^{iv\omega_{k}(0)}-\bar{u}_{k,k}^{iv(0)}|  \quad \forall i\in \mathcal{N}, v\in \mathcal{V}, \omega_{k} \in \Omega_k , \label{m48}\\
& \lambda_{k}^{iv\omega_k(r)} = \rho^{(r-1)}_k|u_{k,k}^{iv\omega_{k}(r)}-\bar{u}_{k,k}^{iv(r-1)}| + \lambda_{k}^{iv\omega_k(r-1)}\quad \forall i\in \mathcal{N}, v\in \mathcal{V}, \omega_{k} \in \Omega_k ,\label{m49} \\
& \sum_{\omega_{k}\in\Omega_k} p(\omega_{k}) (\sum_{i\in \mathcal{N}}\sum_{v\in \mathcal{V}}|u_{k,k}^{iv\omega_{k}}-\bar{u}_{k,k}^{iv}|) \leq \epsilon.\label{m50}
\end{align}

The PHA flow is summarized in Algorithm \ref{alg:PHA}. It iteratively solves penalized subproblems and updates penalty parameters, leading the subproblems under different scenarios to convergence. At the end of every iteration, it will check whether the termination criteria are satisfied. If so, the algorithm will stop, and we can obtain a consensus solution among all scenarios. Note that we also apply the Myopic Policy. Thus, if the routing decision of $\textbf{z}$ reached consensus during the iterations, we will fix the corresponding routing decisions.

\begin{algorithm}[!htb]
\DontPrintSemicolon
Initialize: $r\gets 0, \rho^{(0)}_k$.\\
Solve subproblems \textit{OS}($\omega_{k}$).\\
Calculate consensus variables $\bar{u}_{k,k}^{iv}$ and parameter $\lambda_{k}^{iv\omega_k(0)}$. \\ 
\While{not terminate} {
$r\gets r+1$.\\
Solve subproblems with Objective Function (\ref{m41}).\\
\If{$\textbf{z}$ is consistent}{
Determine the current routing decision as optimal. // Myopic Policy
}
Update consensus variables $\bar{u}_{k,k}^{iv}$, and parameters $\rho^{(r)}_k$ and $\lambda_{k}^{iv\omega_k(r)}$. \\  
}
\Return{$\bar{u}_{k,k}^{iv}$}\\
\caption{The algorithm flow of PHA.}
\label{alg:PHA}
\end{algorithm}

\subsection{Matheuristic Algorithm} \label{5.2}
We adopt the matheuristic algorithm proposed by \citet{solyali2022effective} to solve \textit{OS}($\omega$) and \textit{SS}. It is a heuristic algorithm designed for the Deterministic Inventory Routing Problem (DIRP), mainly based on the problem formulation. The DIRP is first formulated as a Mixed Integer Linear Program (MILP), and then three different MILPs of a restricted version of the formulation are constructed. The matheuristic algorithm relies on a sequential solution of the three MILPs, each continuously optimizing the visiting and delivery plan based on the visiting routes constructed through benchmark algorithms for either the Traveling Salesman Problem (TSP) or the Capacitated Vehicle Routing Problem (CVRP). The main idea of this algorithm is avoiding the subtour elimination constraints, which leads to the main solving difficulty. The specific formulations and implementation within PHA are introduced in Appendix A in our online appendix. We refer the reader to \citet{solyali2022effective} for more details.

\subsection{Heuristic Clustering for Large-Scale Problems} \label{5.3}
Due to the complexity of the data-driven online IRP, even if the scenario decomposition and the matheuristic algorithms are adopted, dealing with large-scale real-life problems in the scenario-optimize step can still be time-consuming. In this subsection, we propose a heuristic clustering algorithm inspired by k-means clustering. Here, we partition retailers into groups with a maximum of $N^{\max}$ retailers at each decision epoch. 

In the IRP decision-making process, the transportation plan is determined not only by considering the relative locations of retailers and the warehouse but also the transportation quantity for each retailer. Thus, our clustering algorithm considers the predicted demands besides the distances between retailers. The clustering criterion, defined as the cost of clustering retailer $i\in\mathcal{N}$ into group $G$, is calculated as \eqref{m111}, where $d(i, G)$ and $d(0, G)$ refer to the geographical distance of retailer $i$ and warehouse $0$ to the center of group $G$, respectively. This equation is basically equivalent to approximating travels with a star centered in the center of group $G$. $r_i$ estimates the utilization of a vehicle of serving retailer $i$, as calculated in \eqref{m222}, where $\epsilon$ is the decay parameter. The delivery quantity is estimated as the sum of decayed forecasted demand (0.5-quantile) for the future $\ell$ periods minus the initial inventory. Then, the fractional utilization bounded between 0 and 1 can be calculated as the fraction of delivery quantity and vehicle capacity. $\lceil \sum_{j\in G}r_j + r_i \rceil -\lceil \sum_{j\in G}r_j \rceil $ indicates whether an additional vehicle must be used if retailer $i$ will join group $G$.

\begin{align}
&c(i, G) = d(i, G) + d(0, G)(\lceil \sum_{j\in G}r_j + r_i \rceil -\lceil \sum_{j\in G }r_j \rceil ), \label{m111}\\
&r_i = \min\{\max\{(\sum_{t=0}^{\ell-1}\epsilon ^t y_{k,k+t}^{\prime i 0.5} - I_{k,k}^i)/Q, 0\}, 1\}.\label{m222}
\end{align}

We apply a standard k-means procedure based on the clustering criteria to partition our retailers into different groups, see for example \cite{sonntag2023stochastic}. We stress that this heuristic clustering is performed based on the results of the scenario-predict step (forecasted 0.5-quantile). Thus, the groups will update at each decision epoch instead of determining a fixed group of retailers for the complete planning horizon.

\section{Computational Experiments} \label{6}
In this section, we first introduce the configuration of our computational experiments in Section 6.1. We create synthetic data, which allows us to control the ``degree of nonstationarity", i.e., meaning that we can manually impose trends and exogenous disturbances in the data. We conduct numerical experiments to show the performance of our ScenPO in Sections 6.2 and 6.3. We end this section by showing the performance of ScenPO on (large-scale) real-life data at our industry partner in Section 6.4. A sensitivity analysis is provided in Appendix B in our online appendix.

\subsection{Experimental Configuration} \label{6.1}
Our ScenPO is implemented in Python, using SciPy \citep{2020SciPy-NMeth}, Scikit-learn \citep{pedregosa2011scikit}, Torch \citep{torch}, and PyTorch-Lightning \citep{falcon2019pytorch}. The MILPs of the matheuristic algorithm are solved using CPLEX version 22.1.1.0 \citep{cplex}. In several parts, for both the ScenPO and benchmark policies introduced later, we employ a Python TSP Solver to solve the TSP \citep{tsp} and a Python wrapper for the Hybrid Genetic Search algorithm proposed by \cite{vidal2012hybrid} to solve capacitated vehicle routing problems \citep{cvrp}. The experiments are run on a computer cluster with 2.4 GHz AMD Genoa 9654 (2x) CPU sockets with 2 GB RAM per core.

In the following, we first introduce our data sets. Afterward, we will introduce seven benchmark types to structurally assess all components of ScenPO. We finally discuss our parameter settings.

\subsubsection{Synthetic data.} We create synthetic instances based on three data patterns: \textit{Trend}, \textit{Random}, and \textit{Both}. For each data pattern, we generate 600 one-year sequences mimicking the demand sequences of 600 retailers, where each sequence is generated as follows. The \textit{Trend} pattern considers a trend and a cycle. We adopt an AR(1) process $v_k = \varphi v_{k-1}$, where for each sequence, $\varphi$ is randomly sampled from a uniform distribution $U(0.995, 1.01)$ and $v_0$ is randomly sampled from the uniform distribution $U(5, 7)$. For the cycle factor, we generate value as $w_k = \sin(2\pi/(\hat{a}_k^y/366)) + \sin(2\pi/(\hat{a}_k^w/7))$, where $\hat{a}_k^y$ and $\hat{a}_k^w$ refer to the day of year and day of week, respectively. The final sequence is generated as $v_k+w_k$. The values of the \textit{Random} pattern are randomly sampled from a normal distribution $N(\mu, \sigma^2)$, where for each sequence, $\mu$ is sampled from $U(9, 11)$ and $\sigma^2$ from $U(2, 4)$. The values of the \textit{Both} pattern are the sum of values generated by the \textit{Trend} and \textit{Random}. We divide the 600 retailers into a training set and a testing set in a ratio of 7:3. The 420 sequences in the training set are fed into prediction models for training, and the testing set is for constructing instances. For the remaining 180 retailers in the testing set, we use k-means clustering to cluster retailers into 10 different groups based on their locations, generated within the latitude of 50.7-53.4 and longitude 3.5-7.1. Finally, we adjust the groups so that each has 18 retailers. Each group corresponds to an instance. Thus, we have 10 instances of each data pattern, each with 18 retailers.

\subsubsection{Benchmarks and ScenPO variants.}
\begin{figure*}[!htb]
    \centering
    \subfigure[MQRNN]{\includegraphics[width=0.32\linewidth]{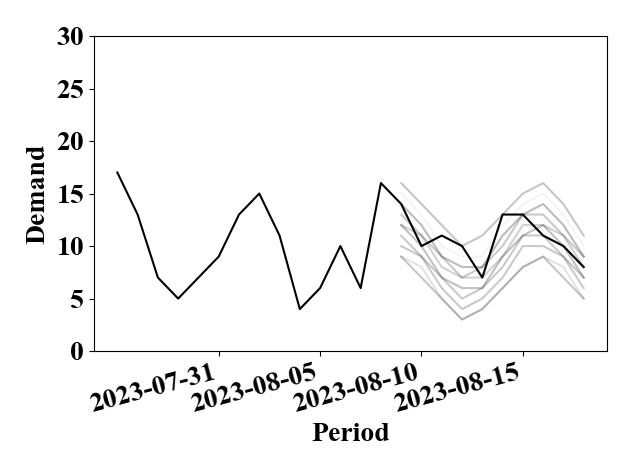}}
    \subfigure[LSTM]{\includegraphics[width=0.32\linewidth]{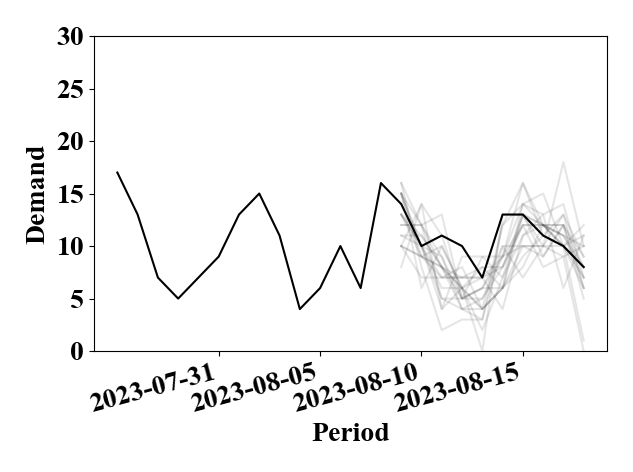}}
    \subfigure[MLE]{\includegraphics[width=0.32\linewidth]{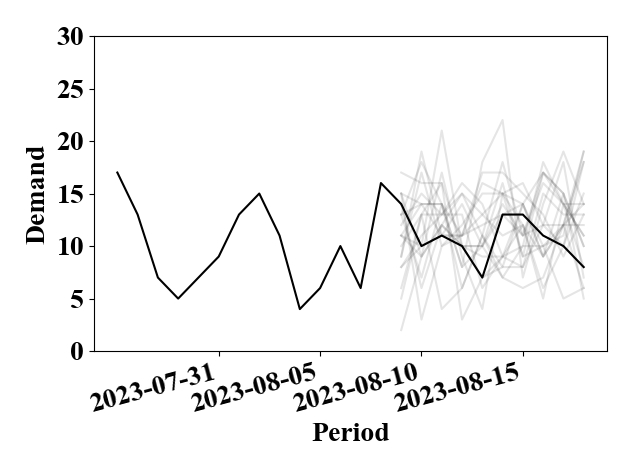}}
    \caption{Demand Scenarios generated via different prediction models on the synthetic data of \textit{Both} data pattern.}
    \label{fig_scenarios}
\end{figure*}

We consider various benchmarks for solving the data-driven online inventory routing problem. These are:
\begin{itemize}
    \item First, we consider a \textit{Perfect Information (PI)} solution, which solves at each epoch $k$ an $\ell$-period single-scenario inventory routing problem, where the single scenario resembles the true future demands. We will compare each other benchmark as the percentage of the gap closed between an expected value solution and the PI solution.
    \item Second, we consider the \textit{Expected Value (EV)} solution, where we also solve at each epoch an $\ell$-period single-scenario inventory routing problem, but we let the demands of the single scenario be equal to the mean demand at each retailer amongst the last $L$ periods.
    \item Third, we consider an \textit{Empirical Sampling (EMP)} solution, in which we randomly sample data of the last $L$ periods to create at each epoch a single-scenario inventory routing problem.
    \item Fourth, we consider a \textit{classic prediction-focused learning paradigm}. Three variants are included, namely where we predict by means of the \textit{MQRNN}, an \textit{LSTM}, or a \textit{Maximum Likelihood Estimation (MLE)}. The single scenario predicted by MQRNN is the 0.5-quantile forecast, $\mathcal{Y}^{\prime 0.5}_{k:}$. We train the LSTM to obtain the point forecast of future demands, $\mathcal{Y}^{\prime LSTM}_{k:}$. For the MLE prediction model, we perform a Chi-square goodness of fit test for normal distribution, exponential distribution, log-normal distribution, Pearson type III distribution, and Weibull minimum distribution on demand sequences to determine the optimal demand distribution and then fit the parameters of the optimal distribution by MLE estimation. At last, the fitted distribution is randomly sampled to generate a single scenario. Summarizing, we consider three prediction-focused learning approaches that predict a single future demand scenario using MQRNN, LSTM, and MLE.
\item Fifth, we employ our \textit{Scenario Predict-then-Optimize paradigm with our SS stochastic programming formulation (ScenPO-SS)}. We consider multiple variants of scenario generation, both quantile- and residuals-based. First, we consider scenario generation using the \textit{MLE}, \textit{LSTM}, and \textit{MQRNN} approaches, as they are also used in the prediction-focused learning paradigm introduced above. For LSTM generating multiple scenarios, we use the point forecast of future demands, $\mathcal{Y}^{\prime LSTM}_{k:}$, as the base value. Then, we sample errors of $\mathcal{E}_{k:}^{\omega_k}$ from a fitted error distribution, which is normal. Demand scenarios are constructed by adding sampled errors to the base value, where different sampled $\mathcal{E}_{k:}^{\omega_k}$ construct different scenarios. For the MLE prediction model, the fitted distribution is randomly sampled multiple times to generate multiple scenarios. The demand scenarios generated by three different models are illustrated in Figure \ref{fig_scenarios}. Note, ScenPO-SS using MQRNN is ``our approach" introduced in this paper. In addition, we consider  the state-of-the-art data-driven Sample Average Approximation (SAA) proposed by \citet{kannan2022data}. They adopt Ordinary Least Squares (OLS) regression to train the prediction model, first used to generate a point prediction. Then, the residuals obtained during training are scaled and added to that point prediction to construct scenarios. They considered three variants of the method, named Empirical Residuals-based SAA (ER-SAA), Jackknife-based SAA (J-SAA), and Jackknife+-based SAA (JP-SAA). We refer readers to \citet{kannan2022data} for more information. We adopt the three approaches as scenario generation benchmarks denoted as \textit{OLS-ER, OLS-J}, and \textit{OLS-JP,} respectively. In addition, we present the performance of the residuals-based approach when we obtain a point prediction using \textit{Holt Winters Exponential Smoothing (called ES)} or using the \textit{0.5-quantile of the MQRNN (called NN)}. It is worth noting that since the ES approach generates predictions by smoothing historical data, it does not require training prediction models to generate point predictions. Thus, we use training data to only obtain residuals by comparing the ground truth to the forecasted demands. For quantile-based methods, in addition to our ScenPO, we adopt \textit{Quantile Regression (QuantReg)} to generate quantiles to evaluate the necessity of deep learning predictors. 
\item Sixth, we consider our \textit{ScenPO approach using the \textit{TS} stochastic programming formulation (ScenPO-TS)}, considering MQRNN, LSTM, and MLE as scenario-generation strategies.
\item Seventh, we consider the state-of-the-art \textit{decision-focused learning approach proposed by \citet{paulus2021comboptnet}, named CombOptNet}. It integrates integer programming solvers into neural network architectures as layers capable of learning both the cost terms and the constraints.
\end{itemize}

\subsubsection{Parameter values.} The per unit holding cost $h$ and per unit backorder cost $e$ are 0.3 and 3, respectively. The scaling parameter $\alpha$ of transportation cost is 0.05 for the synthetic data and 0.01 for the real-life case studies. We use 20 scenarios in our ScenPO variants. The history demand observation of $L$ is set to 14. The look-ahead parameter $\ell$ is set as 5, 7, and 10, respectively. Since the MQRNN directly outputs the $\ell$ period prediction, we need to train separately for the different values of $\ell$. Thus, we train three MQRNN predictors for $\ell=5$, $\ell=7$, and $\ell=10$. LSTM adopts a recursive strategy to generate multi-horizon forecasts, i.e., we train one model to predict a one-step-ahead estimate, then iteratively feed this estimate back as the ground truth to forecast longer periods.

We consider a finite time horizon of 30 periods (or epochs) to evaluate the performance of all the methods. The quantiles forecasted by MQRNN and QuantReg are set as $\mathcal{B}=\{0.1, 0.2, \dots, 0.9\}$. The PHA parameters are set as $\beta_D=1.05, \beta_P=1.05, \rho^{(0)}=0.001$, and $\epsilon=0.1$. As we evaluate over 30 decision epochs, we set the discount factor to one for simplicity. We set a time limit of 900 seconds for solving the MILPs of the PI, EV, EMP, and prediction-focused learning approaches. For the \textit{SS} formulation, we consider a time limit per solved MILP of 1800 seconds, considering the complexity of the formulation. For the \textit{TS} formulation, we set a small time limit of 60 seconds for each MILP in each iteration because we employ PHA to solve it, which requires many more MILP evaluations. We also adopt aggressive PHA convergence parameters to obtain fair comparisons between \textit{SS} and \textit{TS} computationally. All approaches run on the server with a single core, except ScenPO-TS, which runs on 10 cores so that we can solve the subproblems in parallel using the threading module \citep{python-threading}. We set the number of threads invoked by  CPLEX equal to one.  

\subsection{Benchmarking ScenPO}\label{6.2}
In this subsection, we provide a comparison against (a subset of) the seven alternative approaches for solving the data-driven online inventory routing problem. We compare all approaches except the impact of quantile versus residuals-based scenario generation within ScenPO-SS. Those results are presented separately in Section 6.3. The results of CombOptNet are presented in Table A10 in Appendix C in our online appendix as this benchmark approach can only solve small instances due to computational limitations.

\begin{figure}[!tb]
\centering
  \includegraphics[width=0.99\linewidth]{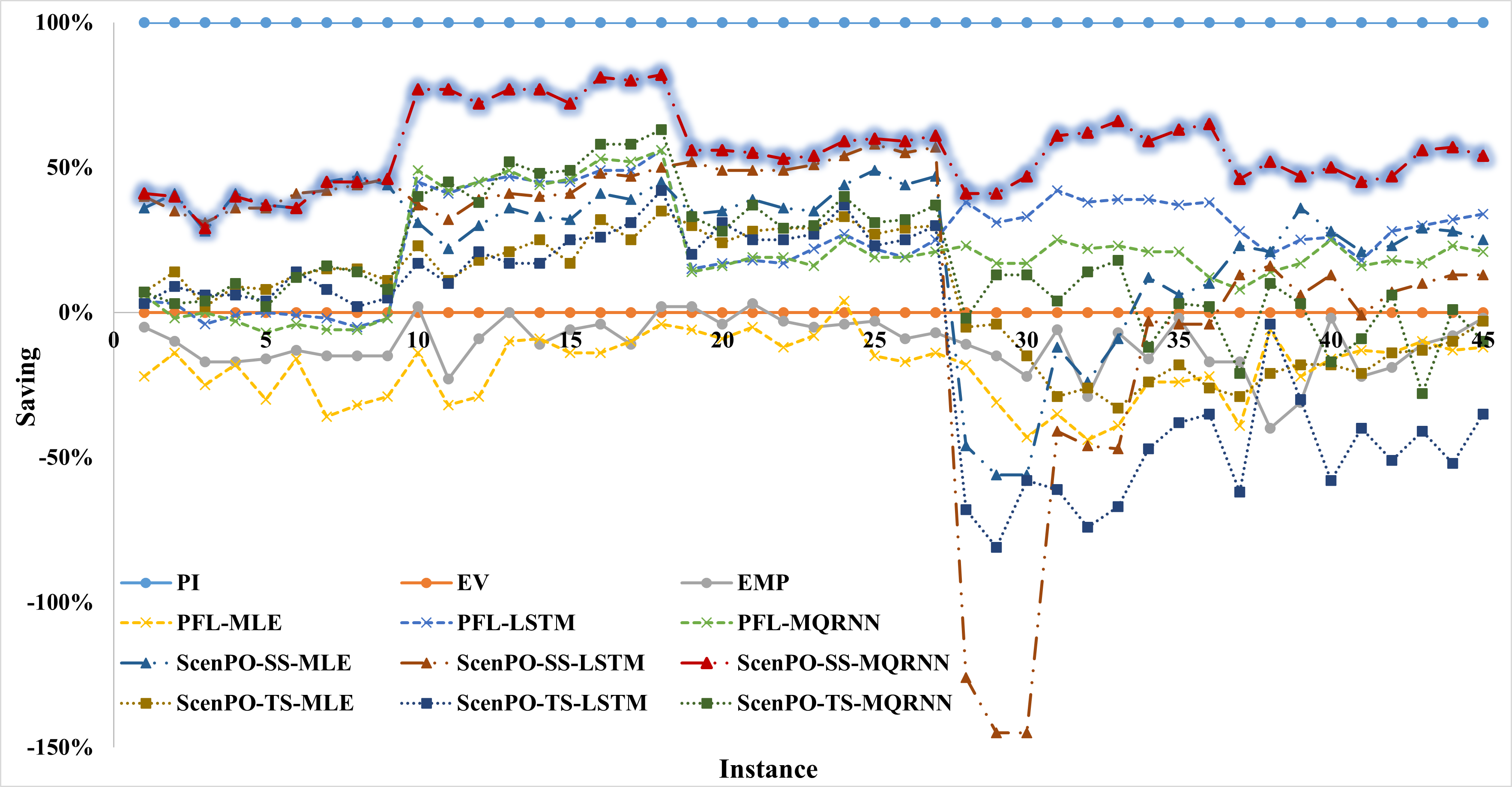}
  \caption{Graphical visualization of all experimental results. The naming of each method refers to the name used in Table \ref{tab:gd}.}
  
  \label{fig_all_exp}
\end{figure}

We present average performances over five instances, and we vary the number of retailers, the look-ahead parameter $\ell$, and the data patterns. A graphical visualization of all the instances, of which details are provided in Tables A4-A6, A11, and A12 in Appendix C in the online appendix, is provided in Figure \ref{fig_all_exp}. Here, we present solely the savings of all methods and arrange all the instances horizontally. Savings means the percentage of the gap between the EV and PI approach being closed by that specific method. The instances are numbered as follows: Instances 1-9 correspond to the \textit{Random} data pattern, Instances 10-18 correspond to the \textit{Trend} data pattern, Instances 19-27 correspond to the \textit{Both} data pattern. The remaining are real-life instances Inst-1 and Inst-2, of which we provide a detailed performance separately in Section 6.4.  Within each data pattern, the instances are numbered as the increase of the number of retailers and $\ell$. For example, Instance 1, Instance 2, and Instance 3 consider five retailers and $\ell=5$, $\ell=7$, and $\ell=10$, respectively. Then Instance 4, Instance 5, and Instance 6 consider seven retailers with $\ell=5$, $\ell=7$, and $\ell=10$, etc. Through the figure, we consistently observe that ScenPO-SS using MQRNN outperforms all other methods. The methods PI, EV, EMP, and the prediction-focused learning methods consider only a single scenario when making decisions. These methods are the fastest to solve, with the average computation time per decision epoch being 705.44s. The ScenPO-SS method takes more time than the single scenario methods, with an average computation time of 1915.11s. The ScenPO-TS method solves formulation \textit{TS} using our PHA algorithm has a computation time of on average 5550.11s.

To further study the detailed differences between the performance of all approaches, we present the results of the instance with ten retailers and the look-ahead length $\ell = 10$ in Table \ref{tab:gd}, for each of the three data patterns. In Table \ref{tab:gd}, the row ``Cost" refers to the total cost over the 30 periods of each approach. The row ``$\Delta$Cost" is the extra cost compared to the PI approach, and ``Saving" means the percentage of the gap between the EV and PI approach being closed. The row ``MSE" is the mean-square error of the used prediction. The row ``Service"  indicates a service level, calculated as the number of directly satisfied demands among all retailers and periods divided by the total demands among all retailers and periods.  Finally, the row ``Time (s)" is the average computational time for each period in seconds.

\begin{table}[hpt]
\centering
\caption{Detailed numerical results of ScenPO and the proposed benchmarks for a 10-retailer system.}
\label{tab:gd}
\scriptsize
\begin{tabular}{c c  r r r r r r r r r r r r}
\hline
\centering
&    & \multirow{2}*{PI} & \multirow{2}*{EV} & \multirow{2}*{EMP} & \multicolumn{3}{c}{Prediction-Focused Learning}         & \multicolumn{3}{c}{ScenPO-SS} &  \multicolumn{3}{c}{ScenPO-TS} \\
\cmidrule(r){6-8}
\cmidrule(r){9-11}
\cmidrule(r){12-14}
&              &            &           &              & MLE         & LSTM    & MQRNN   & MLE     & LSTM    & MQRNN   & MLE              & LSTM     & MQRNN   \\
\hline
\multirow{6}*{\textit{Random}} & Cost    & 1287  & 2339 & 2502  & 2646  & 2360 & 2356 & 1871 & 1852 & 1853 & 2223 & 2289 & 2254 \\
       & $\Delta$Cost   & 0     & 1053 & 1216  & 1360  & 1073 & 1069 & 585  & 566  & 567  & 937  & 1003 & 967  \\
       & Saving  & 100\%     & 0\%    & -15\% & -29\% & -2\% & -2\% & 44\% & \textbf{46\%} & \textbf{46\%} & 11\% & 5\%  & 8\%  \\
       & MSE     & 0     & 275  & 570   & 622   & 290  & 284  &      &      &      &      &      &      \\
       & Service & 100\% & 70\% & 72\%  & 70\%  & 68\% & 68\% & 90\% & 89\% & 88\% & 74\% & 68\% & 72\% \\
       & Time (s)  & 805   & 813  & 810   & 161   & 825  & 814  & 1961 & 1965 & 1741 & 6482 & 6382 & 5666 \\
     \hline
\multirow{6}*{\textit{Trend}}  & Cost    & 1174  & 2087 & 2071  & 2126  & 1579 & 1579 & 1673 & 1628 & 1339 & 1765 & 1705 & 1513 \\
      & $\Delta$Cost   & 0     & 914  & 897   & 952   & 405  & 405  & 499  & 454  & 166  & 591  & 531  & 339  \\
       & Saving  & 100\%     & 0\%    & 2\%   & -4\%  & 56\% & 56\% & 45\% & 50\% & \textbf{82\%} & 35\% & 42\% & 63\% \\
       & MSE     & 0     & 168  & 377   & 466   & 19   & 20   &      &      &      &      &      &      \\
       & Service & 100\% & 70\% & 77\%  & 78\%  & 66\% & 65\% & 95\% & 62\% & 84\% & 84\% & 65\% & 70\% \\
       & Time (s)  & 741   & 762  & 784   & 158   & 762  & 756  & 1925 & 1945 & 1689 & 6085 & 4214 & 3822 \\
     \hline
\multirow{6}*{\textit{Both}}  & Cost    & 1396  & 2606 & 2689  & 2778  & 2298 & 2350 & 2036 & 1918 & 1873 & 2240 & 2245 & 2164 \\
       & $\Delta$Cost   & 0     & 1210 & 1293  & 1382  & 902  & 954  & 640  & 522  & 477  & 844  & 849  & 768  \\
   & Saving  & 100\%     & 0\%    & -7\%  & -14\% & 25\% & 21\% & 47\% & 57\% & \textbf{61\%} & 30\% & 30\% & 37\% \\
       & MSE     & 0     & 319  & 738   & 815   & 177  & 194  &      &      &      &      &      &      \\
       & Service & 100\% & 67\% & 74\%  & 73\%  & 70\% & 67\% & 94\% & 87\% & 88\% & 83\% & 70\% & 74\% \\
       & Time (s)  & 834   & 926  & 863   & 170   & 856  & 858  & 2072 & 2059 & 1879 & 6104 & 5901 & 5295\\
\hline
\end{tabular}
\end{table}

Several observations stand out. First, among the three prediction models, MQRNN performs best, and MLE performs worst. The MLE prediction has the lowest saving, with a -15.67\% saving on average in prediction-focused learning (implying a worse performance than the EV solution), 45.33\% in ScenPO-SS, and 25.33\% in ScenPO-TS. The reason is that this model follows the assumption that the distribution is the same in each period and does not consider the trend and cycle changes in demand. That reveals the shortcomings of straightforward distributional assumptions. The LSTM and MQRNN predictions perform similarly in the prediction-focused learning approach, with average savings of 26.33\% and 25.00\%, respectively. In ScenPO-SS, due to the difference in the scenario generation methods, the MQRNN prediction outperforms the LSTM prediction with a 63.00\% saving, while the average saving of LSTM is 51.00\%. In ScenPO-TS, the average savings of MQRNN and LSTM are 36.00\% and 25.67\%, respectively. The relatively good performance of MQRNN lies in generating scenarios through forecasted quantiles, which is more robust, while LSTM generates scenarios through base values and error distributions.

The second main observation is that ScenPO-SS outperforms all other benchmarks, including ScenPO-TS, with an average saving of 63.00\%. Since MQRNN performs best, we only consider this prediction model in the following analysis. Compared to the EV benchmark, EMP has an average saving of -6.67\%. The saving of ScenPO-SS is 63.00\% on average. In contrast, the average saving of prediction-focused learning, which represents the traditional prediction-focused learning data-driven approach, is only 25.00\%. This shows that ScenPO is a robust alternative that efficiently handles nonstationary data. Relatively surprising, ScenPO-TS does not outperform ScenPO-SS. We will analyze this in more detail later.

Thirdly, ScenPO has the largest improvement under the \textit{Trend} pattern and lowest under the \textit{Random} pattern. This is directly explained by the fact that the \textit{Trend} pattern is the most straightforward to predict. What is also interesting is that the prediction models perform relatively differently among the different data patterns. Focussing on the ScenPO-SS approach, we observe that in the \textit{Trend} pattern, MQRNN outperforms LSTM, having an 82\% saving, and LSTM outperforms MLE with a 50\% saving compared to a 45\% saving. In the \textit{Both} pattern, a similar performance is visible, where the savings of MQRNN, LSTM, and MLE are 61\%, 57\%, and 47\%, respectively. However, the gaps are slightly reduced compared to the \textit{Trend} pattern. In the \textit{Random} pattern, the performance of all models reduces to be similar, with savings of 44\%, 46\%, and 46\% for MLE, LSTM, and MQRNN, respectively. 
\begin{figure}[!htb]
\centering
  \includegraphics[width=0.8\linewidth]{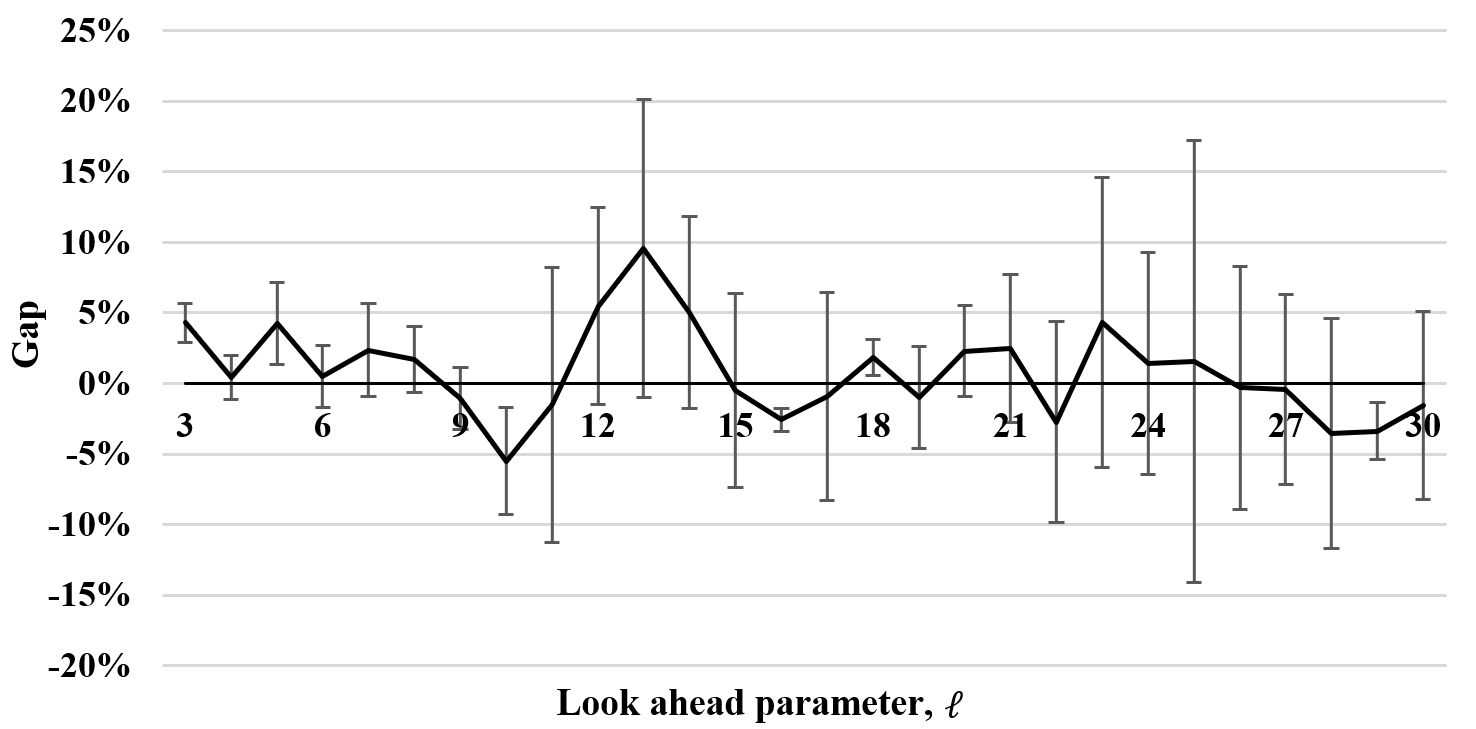}
  \caption{The releative cost comparison of ScenPO-SS compared to ScenPO-TS over different look-ahead lengths $\ell$.}
  \label{fig_1retailer}
\end{figure}
Finally, we noticed that ScenPO-TS does not perform well.  One possible reason is the complexity of the model and the need to use PHA; without that, the models are computationally impossible to solve. Because we need many iterations, we set a one-minute time limit and more aggressive convergence parameters to compress the solution time. To further investigate this, we employ a simple experiment considering the ScenPO-TS and ScenPO-SS approaches for a single retailer setting and vary the value of $\ell$. We set a time limit of 1800s. The experiments are on three instances, and the resulting cost over 30 periods within a rolling horizon setting, for varying levels of $\ell$, is given in Figure \ref{fig_1retailer}, where Gap is calculated as $\frac{{\rm Cost}^{{\rm \textit{TS}}}-{\rm Cost}^{{\rm \textit{SS}}}}{{\rm Cost}^{{\rm \textit{SS}}}}$. We see no clear benefit on either side. ScenPO-SS and ScenPO-TS perform similarly. Future research might want to focus on even more tailored heuristic approaches to solve the $\ell$-period stochastic inventory routing problems to see if the ScenPO-TS approach can outperform the ScenPO-SS approach.  

\subsection{Benchmarking Scenario-Generation Approaches within ScenPO}
To further evaluate the value of scenario generation, we conduct numerical experiments considering different scenario generation benchmarks. We compare various residuals- and quantile-based approaches that can be used within ScenPO to create scenarios in the scenario-generate step. For the scenario-optimize step we use the \textit{SS} stochastic programming formulation.

\begin{figure}[!tb]
\centering
  \includegraphics[width=0.99\linewidth]{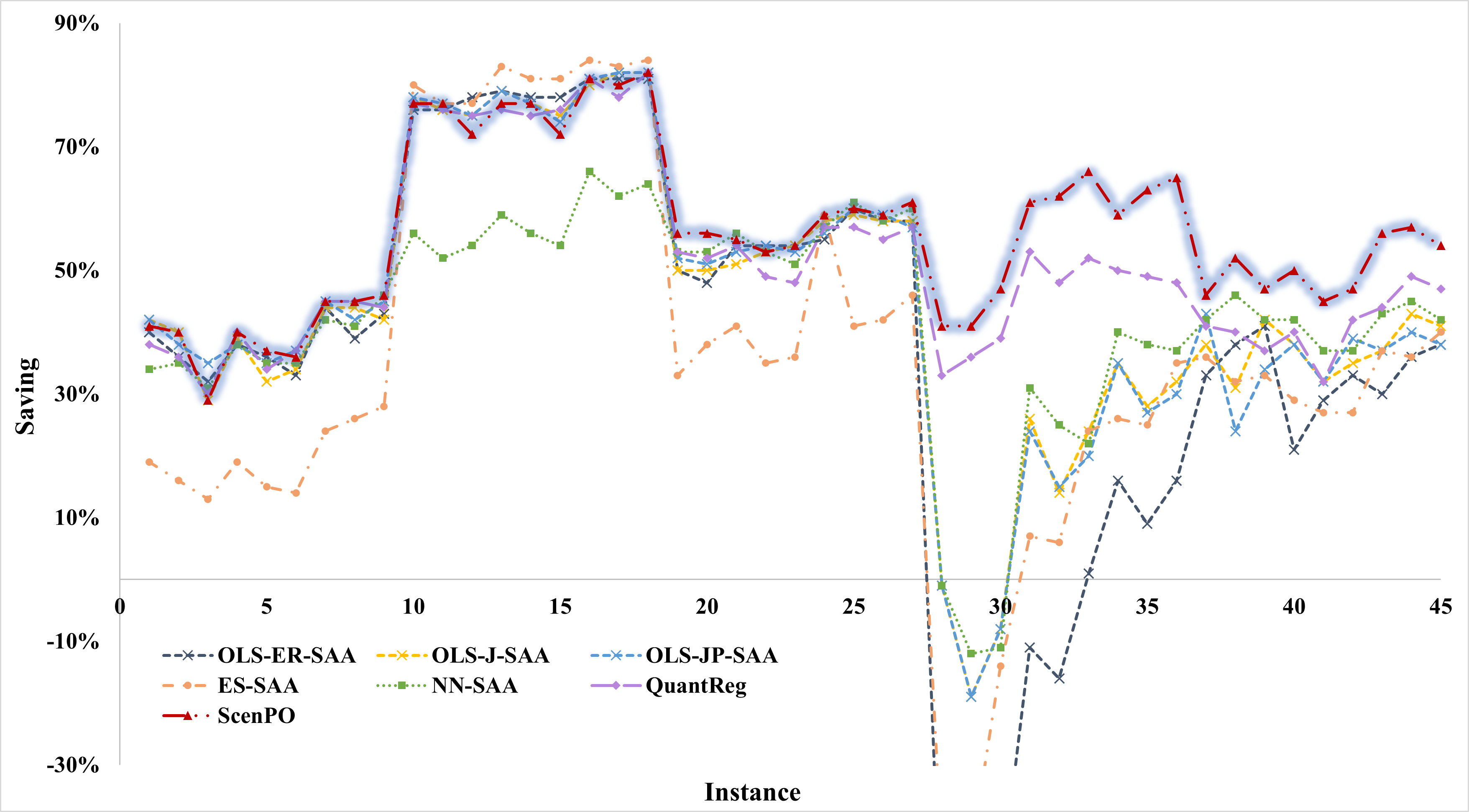}
  \caption{Graphical visualization of the performance of various quantile- and residuals-based scenario generation methods within ScenPO-SS. The naming of each method refers to the name used in Table \ref{tab:expend10_a}.}
  \label{fig_all_exp2}
\end{figure}

We provide a graphical visualization of the performance of all approaches among all the instances in Figure \ref{fig_all_exp2}, similarly as in Figure \ref{fig_all_exp}. We observe that the ``less predictable" the data becomes, the more noticeably the performance benefit of our ScenPO becomes (the blue-shaded red line). Indeed, on the real-life instances (corresponding to Instances 28 and higher in Figure \ref{fig_all_exp2}), our ScenPO approach utilizing MQRNN to generate scenarios outperforms the other scenario-generation variants. We analyze and describe this performance on the real-life instances in Section 6.4.

We further study detailed performances of the residuals-based and quantile-based scenario generation strategies within ScenPO-SS by focusing on the 10-retailer system with $\ell = 10$.
The results are listed in Table \ref{tab:expend10_a}, where the ``ScenPO" column refers to our ScenPO-SS approach with MQRNN as the predictor. Note that all columns, except the PI and EV that serve as a reference, are also ScenPO-SS variants but with different scenario-generation approaches. All other results can be found in Tables A7-A9 in Appendix C in our online appendix.

We obtain the following insights. First, the approaches with simple regression methods as predictors (OLS-ER, OLS-J, OLS-JP, and QuantReg) perform well when the data pattern is clear, that is, under the \textit{Trend} dataset. In such cases, these scenario generation methods perform the same as our ScenPO, which has an 82\% saving. And among the three methods \citet{kannan2022data} proposed, the OLS-ER method performs slightly worse than the other two methods. The exponential smoothing method performs best, achieving 84\% average saving. Second, the different scenario generation methods, except ES, perform similarly when the data pattern has no pattern, that is, under the \textit{Random} dataset. ES performs bad when the data pattern is not clear. Third, when the data pattern is complex, such as under the \textit{Both} dataset, the benchmarks with deep learning-based predictors perform slightly better than simple regression methods. 

Our ScenPO performs well under all datasets and slightly outperforms other benchmarks on our synthetic data in general, with an average of 63.00\%, while the best-performance residuals-based method (OLS-JP) has an average saving of 61.33\%, and QuantReg has an average saving of 61.00\%. As for the computation time, different benchmarks perform similarly, with quantile-based methods taking an average of 1779.5s, slightly faster than residuals-based methods taking 1923.27s on average.

\begin{table}[t]
\centering
\caption{Comparison of the various quantile- and residuals-based scenario generation approaches within ScenPO for a 10-retailer system.}
\scriptsize
\label{tab:expend10_a}
\begin{tabular}{ccrrrrrrrrr}
\hline
  &     & \multirow{2}*{PI} & \multirow{2}*{EV} &  \multicolumn{5}{c}{residuals-based}        & \multicolumn{2}{c}{quantile-based}            \\
\cmidrule(r){5-9}
\cmidrule(r){10-11}
                    &         &            &           & OLS-ER & OLS-J & OLS-JP      & ES &NN & QuantReg          & ScenPO \\
\hline
                     \multirow{5}{*}{\textit{Random}} & Cost         & 1287  & 2339 & 1889            & 1893  & 1867   & 2041 & 1852 & 1872            & 1853   \\
                        &        $\Delta$Cost      & 0     & 1053 & 603             & 606   & 580    & 754  & 566  & 586             & 567    \\
                         &          Saving    & 100\% & 0\%  & 43\%            & 42\%  & 45\%   & 28\% & \textbf{46\%} & 44\%            & \textbf{46\%}   \\
                     &        Service      & 100\% & 70\% & 91\%            & 91\%  & 92\%   & 93\% & 90\% & 89\%            & 88\%   \\
                      &           Time (s)   & 805   & 813  & 1900            & 1920  & 1931   & 1673 & 1906 & 1784            & 1741   \\
                     \hline
                      \multirow{5}{*}{\textit{Trend}}  & Cost         & 1174  & 2087 & 1346            & 1342  & 1340   & 1321 & 1504 & 1336            & 1339   \\
                        & $\Delta$Cost & 0     & 914  & 172             & 168   & 166    & 148  & 330  & 162             & 166    \\
                      & Saving       & 100\% & 0\%  & 81\%            & 82\%  & 82\%   & \textbf{84\%} & 64\% & 82\%            & 82\%   \\
                       & Service      & 100\% & 70\% & 85\%            & 83\%  & 85\%   & 88\% & 72\% & 83\%            & 84\%   \\
                          & Time (s)     & 741   & 762  & 1860            & 1901  & 1902   & 1891 & 1821 & 1691            & 1689   \\
                     \hline
                      \multirow{5}{*}{\textit{Both}}   & Cost         & 1396  & 2606 & 1908            & 1901  & 1916   & 2055 & 1875 & 1917            & 1873   \\
                      & $\Delta$Cost & 0     & 1210 & 512             & 505   & 520    & 659  & 479  & 521             & 477    \\
                      & Saving       & 100\% & 0\%  & 58\%            & 58\%  & 57\%   & 46\% & 60\% & 57\%            & \textbf{61\%}   \\
                      & Service      & 100\% & 67\% & 89\%            & 89\%  & 89\%   & 90\% & 88\% & 86\%            & 88\%   \\
                      & Time (s)     & 834   & 926  & 2084            & 2051  & 2067   & 1845 & 2097 & 1893            & 1879  \\
                     \hline
                     \end{tabular}
\end{table}

\subsection{Real-life Data Experiments}
In this subsection, we apply our approaches to real data concerning aftermarket spare parts demand data from SAIC Volkswagen Automotive Co., Ltd. in China in 2020. The transportation and replenishment for retailers in automotive aftermarket content generally refers to the inventory management and transportation from the original equipment manufacturer (SAIC Volkswagen) to retailers or repair shops (hereafter collectively referred to as retailers). We collect the daily order data from retailers to SAIC Volkswagen to estimate the demands faced by retailers. Currently, SAIC Volkswagen has taken the initiative to establish a real-time inventory information-sharing system with retailers. Our work will be a trusty reference for SAIC Volkswagen to promote their future vendor-manged inventory paradigm. Similar as in \citet{gaur2004periodic}, we aggregate the demands of multiple spare parts into one according to volume. We divide SAIC Volkswagen's 1,000 retailers nationwide into a training set and a testing set according to a 70\%-30\% principle and construct two instances for experiments using a subset of geographically close retailers in the testing set. 

\begin{figure*}[!htb]
    \centering
    \subfigure[\textit{Random}]{\includegraphics[width=0.32\linewidth]{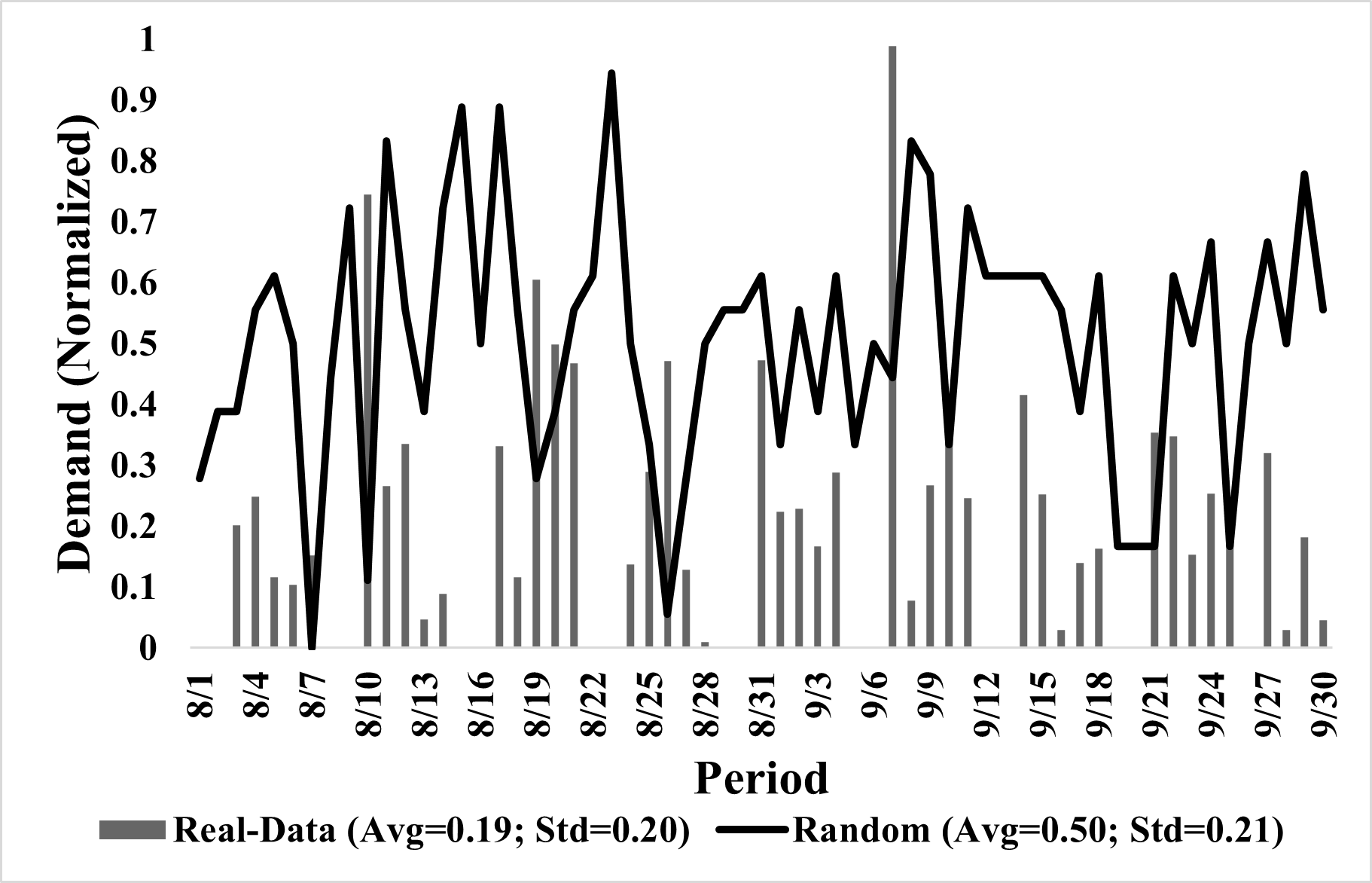}}
    \subfigure[\textit{Trend}]{\includegraphics[width=0.32\linewidth]{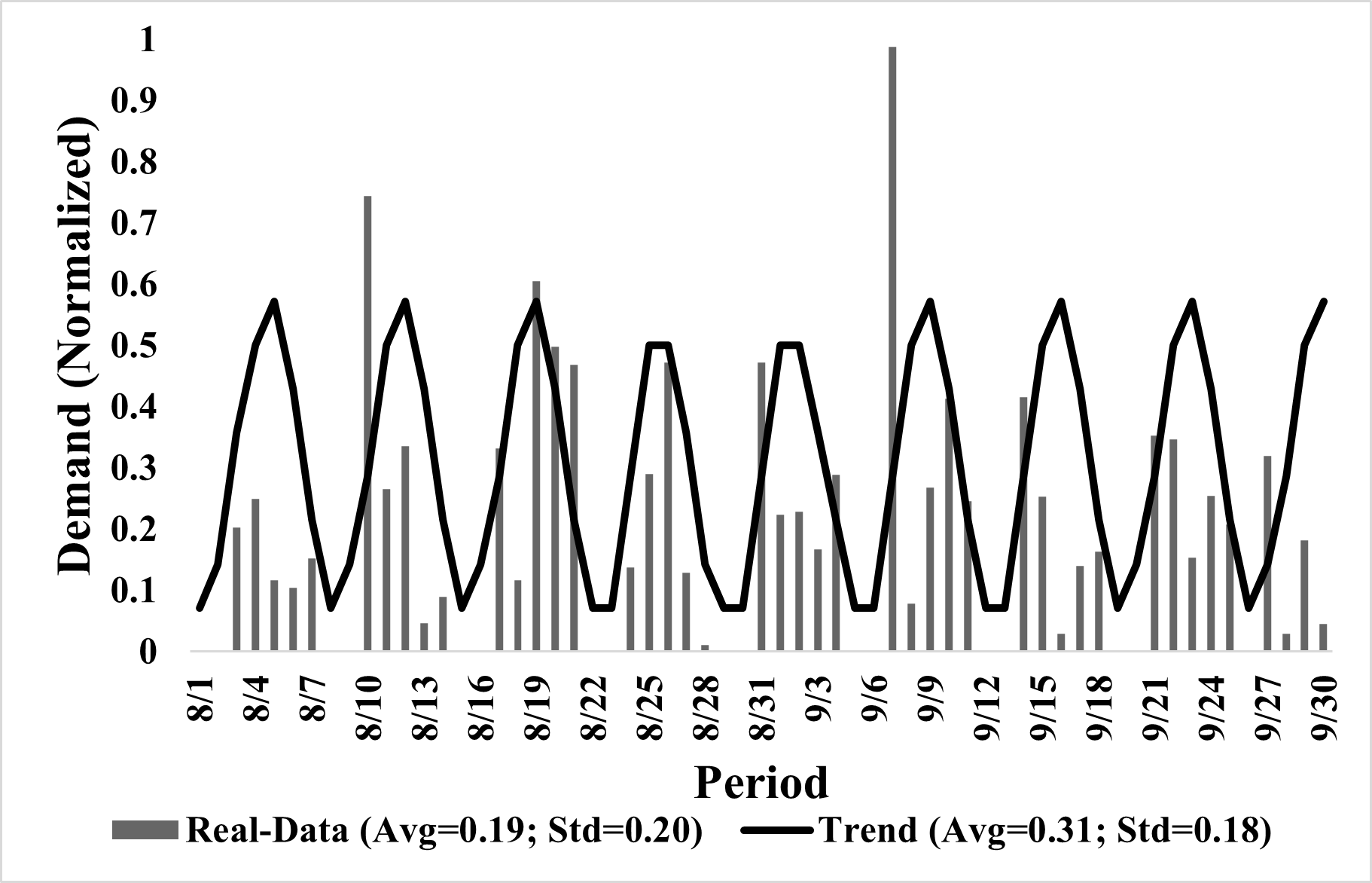}}
    \subfigure[\textit{Both}]{\includegraphics[width=0.32\linewidth]{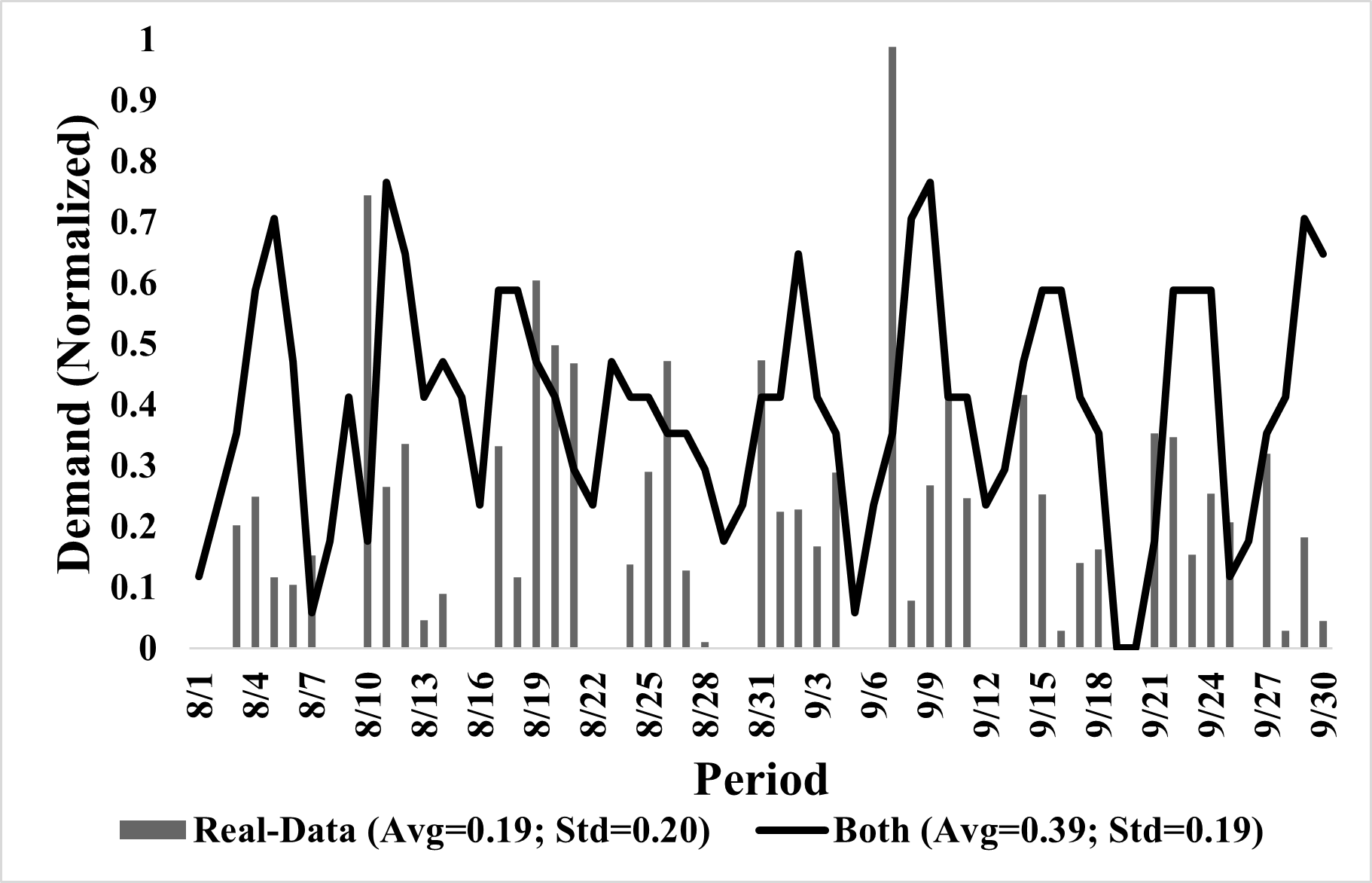}}
    \caption{Comparation of real data with synthetic data. Note, these are not forecasted demand versus actual demands, but solely a visualisation of how real-life demand differs from our synthetic data.}
    \label{fig_real_data}
\end{figure*}
\subsubsection{Example of data used.}
Figure \ref{fig_real_data} shows the normalized demands of one random retailer in August and September, which includes the 30 periods we conduct experiments. In the same plots, we also show the data structure used in our synthetic experiments, being it either \textit{Random}, \textit{Trend}, or \textit{Both}. Compared to our synthetic data,  we observe that the real-life data has a clear pattern that demands usually occur on weekdays, but the specific demand quantity within each week shows no tendency. We can also observe that real-life data fluctuates greatly with some large orders visible.

\begin{table}
\centering
\caption{Real-data experiments for an instance with 10 retailers and $\ell = 10$.}
\label{tab:rd}
\scriptsize
\begin{tabular}{ c  c r r r r r r r r r r r r}
\hline
\centering
& & \multirow{2}*{PI} & \multirow{2}*{EV} & \multirow{2}*{EMP} & \multicolumn{3}{c}{Prediction-Focused Learning} & \multicolumn{3}{c}{ScenPO-SS} &  \multicolumn{3}{c}{ScenPO-TS} \\
\cmidrule(r){6-8}
\cmidrule(r){9-11}
\cmidrule(r){12-14}
&    &            &           &              & MLE         & LSTM    & MQRNN   & MLE & LSTM    & MQRNN   & MLE              & LSTM     & MQRNN   \\
\hline
\multirow{6}*{Inst-1} & Cost    & 585  & 7697  & 8917  & 9296  & 5027  & 6846  & 7006  & 7971  & 3097  & 9571  & 10188 & 7579  \\
       & $\Delta$Cost  & 0    & 7111  & 8332  & 8711  & 4441  & 6261  & 6420  & 7386  & 2511  & 8985  & 9602  & 6993  \\
       & Saving  & 100\%     & 0\%   & -17\% & -22\% & 38\%  & 12\%  & 10\%  & -4\%  & \textbf{65\%}  & -26\% & -35\% & 2\%   \\
       & MSE     & 0    & 137   & 239   & 376   & 50    & 58    &       &       &       &       &       &       \\
       & Service & 98\% & 69\%  & 68\%  & 82\%  & 78\%  & 69\%  & 91\%  & 100\% & 93\%  & 78\%  & 82\%  & 63\%  \\
       & Time (s)  & 87   & 88    & 32    & 19    & 109   & 60    & 50    & 13    & 30    & 5363  & 5201  & 3086  \\
\hline
\multirow{6}*{Inst-2} & Cost    & 6411 & 15103 & 15285 & 16140 & 12130 & 13320 & 12915 & 13981 & 10395 & 15337 & 18126 & 15936 \\
       & $\Delta$Cost  & 0    & 8692  & 8874  & 9729  & 5719  & 6908  & 6504  & 7570  & 3984  & 8926  & 11715 & 9525  \\
       & Saving  &100\%    & 0\%   & -2\%  & -12\% & 34\%  & 21\%  & 25\%  & 13\%  & \textbf{54\%}  & -3\%  & -35\% & -10\% \\
       & MSE     & 0    & 201   & 359   & 410   & 97    & 102   &       &       &       &       &       &       \\
       & Service & 99\% & 68\%  & 74\%  & 74\%  & 77\%  & 72\%  & 89\%  & 99\%  & 89\%  & 79\%  & 79\%  & 66\%  \\
       & Time (s)  & 1194 & 1345  & 633   & 233   & 1348  & 1265  & 912   & 23    & 774   & 10487 & 8931  & 9936\\
\hline
\end{tabular}
\end{table}

\subsubsection{Results on real-life data.} The results of solving the two real-life instances are shown in Tables \ref{tab:rd} and \ref{tab:expandrd1_a}, which consists of 10 retailers with $\ell = 10$. Complete results over multiple combinations of number of retailers and value of $\ell$ can be found in Tables A11-A14 in Appendix C in the online appendix, as well as in the visualization in Figures \ref{fig_all_exp} and \ref{fig_all_exp2}.

Looking at Table \ref{tab:rd}, the different benchmarks perform in the real-life case studies similarly to the synthetic data. MQRNN remains the best under the ScenPO-SS approach, and the average savings of MLE, LSTM, and MQRNN are 17.50\%, 4.50\%, and 59.50\%, respectively. Under the prediction-focused learning approach, the average savings are -17.00\%, 36.00\%, and 16.50\%, with LSTM performing better than MQRNN. ScenPO-SS (with MQRNN) still outperforms all other benchmarks with an average savings of 59.50\%, while prediction-focused learning (with LSTM) has an average saving of 36.00\%, and EMP has an average saving of -9.50\%. We note that the performance of LSTM behaves differently. Under ScenPO-SS, the average saving is 4.50\%, which is worse than MLE (with an average saving of 17.50\%), and the saving is also much lower than LSTM under prediction-focused learning (with an average saving of 36\%). The reason lies in the structure of the real-life industry data: fewer patterns can be captured, and thus, more difficult to predict. That can be evidenced by the gaps of ``Cost" between different benchmarks and PI. The gaps of ScenPO-SS with MQRNN are 43.98\%, 35.28\%, and 40.56\% for \textit{Random}, \textit{Trend}, and \textit{Both} data patterns, respectively, whereas in real-life industry data, the gaps of these two instances are 429.23\% and 62.14\%, respectively. As a result, the error distributions fluctuate significantly, leading to poor quality of the scenarios generated by LSTM. That further illustrates the importance of our MQRNN-based scenario generation method, which is robust enough on different datasets. It also reveals the difficulty of predicting real-life data by point prediction alone. A more powerful multi-scenario approach with robust scenario generation is essential. This is what ScenPO provides.

The computational time to solve the optimization problem in the scenario-optimize step has a similar performance as observed in the synthetic data case, with the approaches considering only a single scenario solving fast, taking on average 534.42s, and the ScenPO-TS method that solves formulation \textit{TS} being slowest with an average computation time of 7167.33s. In contrast, the computation time of ScenPO-SS is 300.33s on average.

\begin{table}[hpt]
\centering
\caption{Further evaluation of different scenario generation approaches for real-date instance with 10 retailers and $\ell = 10$.}
\label{tab:expandrd1_a}
\scriptsize
\begin{tabular}{ccrrrrrrrrr}
\hline
  &     & \multirow{2}*{PI} & \multirow{2}*{EV} & \multicolumn{5}{c}{residuals-based}        & \multicolumn{2}{c}{quantile-based}            \\
\cmidrule(r){5-9}
\cmidrule(r){10-11}
                    &         &            &          & OLS-ER & OLS-J & OLS-JP  &ES    & NN & QuantReg          & ScenPO \\
\hline
         \multirow{5}*{Inst-1}       & Cost         & 585   & 7697    & 6551            & 5455  & 5551   & 5240  & 5091  & 4251            & 3097   \\
                    & $\Delta$Cost & 0     & 7111    & 5966            & 4870  & 4965   & 4654  & 4506  & 3666            & 2511   \\
                    & Saving       & 100\% & 0\%     & 16\%            & 32\%  & 30\%   & 35\%  & 37\%  & 48\%            & \textbf{65\%}   \\
                    & Service      & 98\%  & 69\%    & 93\%            & 93\%  & 92\%   & 93\%  & 94\%  & 88\%            & 93\%   \\
                    & Time (s)     & 87    & 88      & 51              & 70    & 73     & 45    & 106   & 23              & 30    \\
           \hline
 \multirow{5}{*}{Inst-2} & Cost         & 6411  & 15103 & 11805           & 11514 & 11826  & 11667 & 11485 & 11015           & 10395  \\
                     & $\Delta$Cost & 0     & 8692  & 5394            & 5103  & 5415   & 5256  & 5074  & 4604            & 3984   \\
                     & Saving       & 100\% & 0\%   & 38\%            & 41\%  & 38\%   & 40\%  & 42\%  & 47\%            & \textbf{54\%}   \\
                     & Service      & 99\%  & 68\%  & 92\%            & 88\%  & 89\%   & 90\%  & 91\%  & 85\%            & 89\%   \\
                     & Time (s)     & 1194  & 1345  & 1167            & 1864  & 1948   & 1298  & 2414  & 628             & 774   \\
                     \hline
\end{tabular}
\end{table}
                     
The performance of scenario generation benchmarks is evaluated in Table \ref{tab:expandrd1_a} for the 10-retailer system with $\ell = 10$. Complete results over multiple parameter combinations can be found in Tables A13 and A14 in Appendix C in the online appendix. Our ScenPO outperforms the other benchmarks, with an average saving of 60.00\%, while the best-performance residuals-based approach (NN) has an average saving of 40.00\%. This results stand out, as the NN approach takes the 0.5-quantile prediction of the MQRNN, but then samples scenarios via the residuals. It can be interpreted as the residuals-based variant of the method we propose in this paper, as we also use the MQRNN but then sample scenarios uniformly over the quantiles as predicted by the MQRNN.

In addition, we notice that QuantReg also performs competitively in some instances, and the average saving is 48.00\%, illustrating the value of the quantile-based scenario generation approach as it also outperforms NN. We also find that taking the deep-learning method as a predictor performs better than simple regression methods, whether for residuals-based or quantile-based methods, showing that deep-learning methods are advantageous to capture the complex correlations for complex real-life data. As for the computation time, the quantile-based methods still have the advantage, with the average time being 363.75s, while the residuals-based methods take 903.60s on average.

\subsubsection{Performance on large-scale real-life data.} Finally, since real-life applications always concern many retailers,  we also construct two large-scale instances that consider the transportation and replenishment plan for 100 retailers based upon our real-life data from our industry partner. We utilize the heuristic clustering algorithm introduced in Section \ref{5.3}, where the maximum number of retailers in each group $N^{\max}$ is set to seven. The results of the large-scale problems are shown in Table \ref{tab:100r}. The results show that our ScenPO still performs consistently well on large-scale instances. Compared to classic prediction-focused learning with the best-performed predictor (LSTM), our best-performed ScenPO-SS (with predictor MQRNN) decreases the cost by 15.74\% on average, and achieves a service level over 90\%. Compared to the best-performance residuals-based method (NN), our best-performed ScenPO-SS decreases the cost by 15.18\%. QuantReg also performs well, reducing the cost by 8.64\% compared to the best-performance residuals-based method. Regarding computation times, ScenPO-SS and QuantReg perform better than other methods and take 69.92s on average. The residuals-based method takes 539.90s on average, while prediction-focused learning approaches take an average of 449.56s.

\begin{table}
\centering
\scriptsize
\caption{Real-data experiments of 100-retailer instances.}
\label{tab:100r}
\begin{tabular}{c c c r r r r r r r r r r r r}
\hline
\centering
\multirow{2}*{Instance} & \multirow{2}*{$\ell$} & \multirow{2}*{Item}         & \multicolumn{3}{c}{Prediction-Focused Learning}         & \multicolumn{3}{c}{ScenPO-SS} & \multicolumn{5}{c}{residuals-based}  &\multirow{2}*{QuantReg}\\
\cmidrule(r){4-6}
\cmidrule(r){7-9}
\cmidrule(r){10-14}
                            &                                      &                       & MLE             & LSTM          & MQRNN         & MLE     & LSTM     & MQRNN  & OLS-ER & OLS-J & OLS-JP & ES    & NN    &         \\
                            \hline
\multirow{9}{*}{L\_Inst\_1} & \multirow{3}{*}{5}                   & Cost                  & 116894          & 78372         & 85667         & 97368   & 119430   & \textbf{65607}  & 84697  & 78216 & 77494  & 89931 & 75729 & 71096           \\
                            &                                      & Service               & 83\%            & 81\%          & 74\%          & 94\%    & 99\%     & 92\%   & 94\%   & 92\%  & 92\%   & 95\%  & 93\%  & 89\%            \\
                            &                                      & Time (s)              & 17              & 183           & 76            & 30      & 11       & 17     & 57     & 154   & 112    & 35    & 186   & 14              \\
                            \cline{2-15}
                            & \multirow{3}{*}{7}                   & Cost                  & 118115          & 78171         & 87154         & 96012   & 120991   & \textbf{66580}  & 92334  & 82985 & 83378  & 88287 & 78599 & 69893           \\
                            &                                      & Service               & 84\%            & 82\%          & 74\%          & 93\%    & 99\%     & 92\%   & 95\%   & 93\%  & 93\%   & 96\%  & 93\%  & 89\%            \\
                            &                                      & Time (s)              & 221             & 736           & 363           & 97      & 14       & 25     & 169    & 562   & 484    & 257   & 1112  & 25              \\
                            \cline{2-15}
                            & \multirow{3}{*}{10}                  & Cost                  & 112032          & 77575         & 83239         & 96071   & 120708   & \textbf{66339}  & 91275  & 78972 & 77887  & 86240 & 77331 & 71299           \\
                            &                                      & Service               & 85\%            & 82\%          & 75\%          & 94\%    & 99\%     & 91\%   & 94\%   & 92\%  & 92\%   & 94\%  & 93\%  & 88\%            \\
                            &                                      & Time (s)              & 858             & 1213          & 1174          & 184     & 25       & 172    & 962    & 1394  & 1248   & 832   & 2248  & 102             \\
                            \hline
\multirow{9}{*}{L\_Inst\_2} & \multirow{3}{*}{5}                   & Cost                  & 102252          & 67518         & 75885         & 86293   & 108049   & \textbf{56473}  & 74867  & 67826 & 67874  & 78692 & 65914 & 60960           \\
                            &                                      & Service               & 83\%            & 80\%          & 72\%          & 94\%    & 99\%     & 92\%   & 95\%   & 94\%  & 93\%   & 97\%  & 94\%  & 89\%            \\
                            &                                      & Time (s)              & 22              & 47            & 30            & 29      & 10       & 19     & 24     & 27    & 33     & 30    & 101   & 11              \\
                            \cline{2-15}
                            & \multirow{3}{*}{7}                   & Cost                  & 102801          & 67574         & 75731         & 84948   & 109426   & \textbf{57124}  & 83232  & 70862 & 70948  & 77586 & 67681 & 61199           \\
                            &                                      & Service               & 83\%            & 81\%          & 73\%          & 94\%    & 99\%     & 92\%   & 95\%   & 94\%  & 94\%   & 97\%  & 94\%  & 88\%            \\
                            &                                      & Time (s)              & 101             & 322           & 137           & 32      & 14       & 41     & 90     & 225   & 228    & 112   & 441   & 28              \\
                            \cline{2-15}
                            & \multirow{3}{*}{10}                  & Cost                  & 101753          & 67605         & 73341         & 86709   & 109406   & \textbf{56112}  & 81439  & 68655 & 68390  & 74965 & 68681 & 61979           \\
                            &                                      & Service               & 83\%            & 81\%          & 73\%          & 94\%    & 99\%     & 92\%   & 95\%   & 92\%  & 92\%   & 96\%  & 94\%  & 88\%            \\
                            &                                      & Time (s)              & 712             & 989           & 891           & 84      & 19       & 240    & 791    & 1058  & 990    & 505   & 1730  & 145      \\
                            \hline
\end{tabular}
\end{table}

\section{Conclusions} \label{7}
This paper proposes a novel integrated prediction and optimization approach for real-time joint replenishment and distribution via vehicle routing in one-warehouse, multiple-retailer systems, also known as data-driven online inventory routing. Our approach, Scenario Predict-then-Optimize (ScenPO), relies neither on making distributional assumptions for retailer demand nor on single-point forecasts. This contrasts the state-of-the-art that makes such assumptions or works with point forecasts. ScenPO consists of a \textit{scenario-predict} step, utilizing the multi-horizon quantile recurrent neural network that forecasts future retailer demand quantiles from which we sample scenarios. The sampled scenarios are input for the \textit{scenario-optimize} step of ScenPO, which exploits two-stage stochastic programming to find a solution that minimizes the future expected cost of replenishment and distribution decisions.

In this way, ScenPO can capture (auto)correlation and lumpy retailer demand patterns that typically occur in one warehouse and multiple retailers systems. We show our ScenPO's efficiency based on synthetic data and large-scale real-life data from SAIC Volkswagen Automotive Co., Ltd. We compare ScenPO to empirical sampling techniques, various classical prediction-focused learning approaches, and state-of-the-art data-driven approaches including the decision-focused learning approach, CombOptNet, and various residuals-based sample average approximation approaches for generating scenarios. The results show that our ScenPO outperforms these approaches significantly. It closes 61.60\% of the gap between the expected value solution and an oracle solver, whereas the best alternative approach only closes 55.60\% on average.

Our ScenPO is general because it leverages neural network-based quantile prediction and stochastic programming. It can easily be tailored to other applications than online inventory routing in which nonstationary demand patterns and online combinatorial optimization are combined. Future research is needed to investigate how to best align the scenario-predict and scenario-optimize steps in other application areas. 

Another interesting future research direction is to investigate distributionally robust approaches for application in which only limited data is available. The complexity of the inventory routing problem warrants tailored approaches based upon such distributionally robust approaches, specifically in a dynamic environment.

Another promising avenue for further research is to enrich ScenPO for other variants of online inventory routing problems. For example, one might investigate cyclic delivery patterns to enhance the alignment with other planning processes, consider limited supply availability at the central warehouse, or consider more complex supply chains arising in, for instance,  modern city logistics landscapes due to environmental bans of heavy-duty vehicles.
 
\bibliographystyle{informs2014}
\bibliography{reference}

\clearpage
\begin{APPENDICES}  
\setcounter{equation}{0}
\renewcommand{\theequation}{A\arabic{equation}}

\setcounter{table}{0}   
\setcounter{figure}{0}

\renewcommand{\thetable}{A\arabic{table}}
\renewcommand{\thefigure}{A\arabic{figure}}

\section{Matheuristic Algorithm} \label{a1}
The algorithm follows the following steps:  
\begin{itemize}
  \item[Step 1.] Solve a TSP over all retailers,
  \item[Step 2.] Solve \textit{MILP}1,
  \item[Step 3.] Solve a TSP if only one vehicle is enabled or a CVRP else for each $t\in \mathcal{T}$,
  \item[Step 4.] Solve \textit{MILP}2,
  \item[Step 5.] Solve a TSP for each route in each $t\in \mathcal{T}$,
  \item[Step 6.] Solve \textit{MILP}3,
  \item[Step 7.] Solve a TSP for each vehicle in each $t\in \mathcal{T}$.
\end{itemize}
Next, we will introduce the MILPs and the penalized objective functions implemented in the PHA. Here, we omit the notation of $\omega_k$.

Based on the TSP solution, we define $\alpha_k(i)$ and $\beta_k(i)$ as the retailer sets visited after and before retailer $i$. When $i=0$ or $j=0$, $z_{k,k+t}^{ij}$ is an integer variable indicating the visit times to the warehouse. Otherwise, $z_{k,k+t}^{ij}$ is a binary variable indicating whether the corresponding edge is visited. Then, we can formulate the \textit{MILP}1 formulation. Compared to \textit{OS}($\omega_{k}$), we ignore the vehicle index $v$ in the \textit{MILP}1 formulation and construct the vehicle route based on the TSP solution. Thus, the sub-tour elimination constraints are no longer required. 

\begin{align}
(\textit{MILP}1) \min\ & \sum_{t\in \mathcal{T}}\sum_{(i, j)\in E}  c_{ij}z_{k,k+t}^{ij} +\sum_{t\in \mathcal{T}}\sum_{i\in \mathcal{N}}(h (\hat{I}_{k,k+t+1}^{i})^+ + \!\!\!\!\!\!\!\!\!\!\!\!\!\!\!\!\!\!\!\!\!\!\!\!\!\!\!\!\!\!\!\!\!\!\!\!\!\!\!\!\!\!\!\!\!\!\!\!\!\!\!&& e (\hat{I}_{k,k+t+1}^i)^-),\\
  s.t. \ &\sum_{(i,j)\in \xi^+(0)} z_{k,k+t}^{ij}=\sum_{(i,j)\in \xi^-(0)} z_{k,k+t}^{ij}= V\!\!\!\!\!\!\!\!\!\!\!\!\!\!\!\!\!\!\!\!\!\!\!\!\!\!\!\!\!\!\!\!\!\!\!\!\!\!\!\!\!\!\!\!\!\!\!\!\!\!\!\! &&\forall t\in \mathcal{T} ,\\
  &\sum_{j\in \alpha(i)} z_{k,k+t}^{ij} = \sum_{j\in \beta(i)} z_{k,k+t}^{ij} \!\!\!\!\!\!&&\forall i\in \mathcal{N}, t\in \mathcal{T},\\
  &\sum_{j\in \alpha(i)} z_{k,k+t}^{ij} \leq 1 &&\forall i\in \mathcal{N}, t\in \mathcal{T},\\
  &u_{k,k+t}^{i} \leq M\sum_{j\in \alpha(i)} z_{k,k+t}^{ij} &&\forall i\in \mathcal{N}, t\in \mathcal{T} ,\\  
  &\sum_{i\in \mathcal{N}} u_{k,k+t}^{i} \leq Q \sum_{j\in \mathcal{N}} z_{k,k+t}^{0j} \!\!\!\!&&\forall t\in \mathcal{T} ,\\
  &u_{k,k+t}^{i} \leq I^{\max} - \hat{I}_{k,k+t}^i &&\forall i\in \mathcal{N}, t\in \mathcal{T} ,\\
  &\hat{I}_{k,k+t+1}^i = \hat{I}_{k,k+t}^i + u_{k,k+t}^{i} - \hat{y}_{k,k+t}^i \!\!\!\!\!\!\!\!\!\!\!\!\!\!\!\!\!\!\!\!\!\!\!\!\!\!\!\!\!\!\!&&\forall i\in \mathcal{N}, t\in \mathcal{T} ,\\ 
  &z_{k,k+t}^{ij}\in \{0,1\} &&\forall (i,j)\in E, t\in \mathcal{T},\\
  &u_{k,k+t}^{i} \geq 0 &&\forall i\in \mathcal{N}, t\in \mathcal{T}.
\end{align}

When implemented in PHA, the penalized objective will be modified as Equation (\ref{mA12}), while constraints (\ref{mA13}) and (\ref{mA14}) are added. The aggregated PHA parameter $\lambda_{k}^{i}$ is calculated as Equation (\ref{mA15}) and (\ref{mA16}).

\begin{align}
\min\ & \sum_{t\in \mathcal{T}}\sum_{(i, j)\in E}  c_{ij}z_{k,k+t}^{ij} +\sum_{t\in \mathcal{T}} \sum_{i\in \mathcal{N}}(h( \hat{I}_{k,k+t+1}^i)^+ + e (\hat{I}_{k,k+t+1}^i)^-) + \sum_{i\in \mathcal{N}}\lambda_k^{i}(\mu_{k,k}^{i}+\nu^{i}_{k,k}), \label{mA12}\\
s.t. & \ \mu_{k,k}^{i}-\nu^{i}_{k,k} = u_{k,k}^{i}-\sum_{v\in \mathcal{V}}\bar{u}_{k,k}^{iv} \quad \forall i\in \mathcal{N}, t\in \mathcal{T}, \label{mA13}\\
& \mu_{k,k}^{i},\nu^{i}_{k,k} \geq 0  \qquad \qquad \qquad \quad \ \ \forall i\in \mathcal{N}, t\in \mathcal{T}. \label{mA14}
\end{align}

\begin{align}
& \lambda_k^{i(0)}= \rho_k^{(0)}|\sum_{v\in\mathcal{V}}u_{k,k}^{iv(0)}-\sum_{v\in\mathcal{V}}\bar{u}_{k,k}^{iv(0)}| \quad \forall i\in \mathcal{N},  \label{mA15}\\
& \lambda_k^{i(r)} = \rho_k^{(r-1)}|\sum_{v\in\mathcal{V}}u_{k,k}^{iv(r)}-\sum_{v\in\mathcal{V}}\bar{u}_{k,k}^{iv(r-1)}| + \lambda_k^{i(r-1)}  \quad \forall i\in \mathcal{N}.  \label{mA16}
\end{align}

Since \textit{MILP}1 only considers the aggregated vehicle capacity constraints, the solution obtained may be unfeasible. Therefore, \textit{MILP}2 is constructed to make the solution feasible or to improve it by exchanging retailers among routes. The notations needed are defined as follows.

\begin{table}[!htb]
  \centering
  \small
  \caption{Additional notations for \textit{MILP}2.}
  \scalebox{1}{\begin{tabular}{ll}
  \hline
    \multicolumn{2}{l}{\textbf{Parameters}} \\
    $\Upgamma_{k,k+t}^{iv}$ & driving cost of inserting retailer $i$ to route $v$ in $k+t$\\
    $\Updelta_{k,k+t}^{iv}$ & driving cost saving of removing retailer $i$ from route $v$ in $k+t$\\
    $a_{k,k+t}^{iv}$ & indicator parameter of whether retailer $i$ is visited by route $v$ in $k+t$ in the solution\\
    \multicolumn{2}{l}{\textbf{Decision Variables}} \\  
    $s_{k,k+t}^{iv}$ & binary variable indicating whether retailer $i$ is inserted to route $v$ in $k+t$\\
    $r_{k,k+t}^{iv}$ & binary variable indicating whether retailer $i$ is removed from route $v$ in $k+t$\\
  \hline
  \end{tabular}}
  \label{tab: additional notations1}
\end{table}

\begin{align}
(\textit{MILP}2) \min\ &\sum_{t\in \mathcal{T}} \sum_{i\in \mathcal{N}}(h (\hat{I}_{k,k+t+1}^i)^+ + e (\hat{I}_{k,k+t+1}^i)^-) + \sum_{t\in \mathcal{T}} \sum_{i\in \mathcal{N}} \sum_{v\in \mathcal{V}}((1-a_{k,k+t}^{iv})\!\!\!\!\!\!\!&&\Upgamma_{k,k+t}^{iv} s_{k,k+t}^{iv} - a_{k,k+t}^{iv} \Updelta_{k,k+t}^{iv} r_{k,k+t}^{iv}) ,\\
  \ s.t. \ &\sum_{v\in \mathcal{V}} (a_{k,k+t}^{iv} - a_{k,k+t}^{iv} r_{k,k+t}^{iv} + (1-a_{k,k+t}^{iv}) s_{k,k+t}^{iv}) \leq 1 &&\forall i\in \mathcal{N}, t\in \mathcal{T},\\
  &u_{k,k+t}^{iv} \leq M(a_{k,k+t}^{iv} - a_{k,k+t}^{iv} r_{k,k+t}^{iv} + (1-a_{k,k+t}^{iv}) s_{k,k+t}^{iv})  &&\forall i\in \mathcal{N}, t\in \mathcal{T},  v\in \mathcal{V} ,\\
  &\sum_{i\in \mathcal{N}} u_{k,k+t}^{iv} \leq Q &&\forall t\in \mathcal{T}, v\in \mathcal{V} ,\\
  &u_{k,k+t}^{iv} \leq I^{\max} - \hat{I}_{k,k+t}^{i} &&\forall i\in \mathcal{N}, t\in \mathcal{T}, v\in \mathcal{V} ,\\
  &\hat{I}_{k,k+t+1}^i = \hat{I}_{k,k+t}^{i} + \sum_{v\in \mathcal{V}}u_{k,k+t}^{iv} - \hat{y}_{k,k+t}^i &&\forall i\in \mathcal{N}, t\in \mathcal{T} ,\\
  &r_{k,k+t}^{iv} \in \{0, 1\} &&\forall i\in \mathcal{N}, t\in \mathcal{T}, v\in \mathcal{V},\\
  &s_{k,k+t}^{iv} \in \{0, 1\} &&\forall i\in \mathcal{N}, t\in \mathcal{T}, v\in \mathcal{V},\\
  &u_{k,k+t}^{iv} \geq 0 &&\forall i\in \mathcal{N}, t\in \mathcal{T}, v\in \mathcal{V} .
\end{align}

The objective of \textit{MILP}2 in PHA is modified to Equation (\ref{mA25}), with Constraints (\ref{mA26}) and (\ref{mA27}). 

\begin{align}
\min \ & \sum_{t\in \mathcal{T}} \sum_{i\in \mathcal{N}}(h (\hat{I}_{k,k+t+1}^i)^+ + e (\hat{I}_{k,k+t+1}^i)^-) + \sum_{t\in \mathcal{T}} \sum_{i\in \mathcal{N}} \sum_{v\in \mathcal{V}}((1-a_{k,k+t}^{iv})\Upgamma_{k,k+t}^{iv} s_{k,k+t}^{iv} - a_{k,k+t}^{iv} \Updelta_{k,k+t}^{iv} r_{k,k+t}^{iv})+\nonumber\\
& \sum_{i\in \mathcal{N}}\sum_{v\in \mathcal{V}}\lambda_{k}^{iv}(\mu_{k,k}^{iv}+\nu_{k,k}^{iv}), \label{mA25}\\
s.t. \  & \mu_{k,k}^{iv}-\nu_{k,k}^{iv} = u_{k,k+t}^{iv}-\bar{u}_{k,k+t}^{iv} \quad \forall i\in \mathcal{N}, v\in \mathcal{V}, \label{mA26}\\
& \mu_{k,k}^{iv}, \nu_{k,k}^{iv} \geq 0  \qquad \qquad \qquad \quad \ \forall i\in \mathcal{N}, v\in \mathcal{V} .\label{mA27}
\end{align}

Then, formulation \textit{MILP}3 is further formulated to enhance the solution with the notations listed in Table \ref{tab: additional notations2}.

\begin{table}[!htb]
  \centering
  \small
  \caption{Additional notations for \textit{MILP}3.}
  \scalebox{1}{\begin{tabular}{ll}
  \hline
    \multicolumn{2}{l}{\textbf{Sets}} \\
    $\mathcal{V}_{k,k+t}$ & set of vehicles used in $k+t$ in the solution of the last step, $\mathcal{V}_{k,k+t}=\{1, 2, \dots, V_{k,k+t}\}$\\
    \multicolumn{2}{l}{\textbf{Parameters}} \\
    $\alpha_{k,k+t}^{v}(i)$ & the retailer set visited after retailer $i$ in $k+t$ by route $v$\\
    $\beta_{k,k+t}^{v}(i)$ & the retailer set visited before retailer $i$ in $k+t$ by route $v$\\
  \hline
  \end{tabular}}
  \label{tab: additional notations2}
\end{table}

\begin{align}
(\textit{MILP}3) \min\ & \sum_{t\in \mathcal{T}}\sum_{(i,j)\in E}\sum_{v\in \mathcal{V}_{k,k+t}}  c_{ij}z_{k,k+t}^{ijv}+ \sum_{t\in \mathcal{T}} \sum_{i\in \mathcal{N}}(h (\hat{I}_{k,k+t+1}^i)^+ + e \!\!\!\!\!\!\!\!\!\!\!\!&& (\hat{I}_{k,k+t+1}^i)^-), \\
  s.t. \ &\sum_{(i,j)\in \xi^+(0)} z_{k,k+t}^{ijv}=\sum_{(i,j)\in \xi^-(0)} z_{k,k+t}^{ijv} = 1 &&\forall t\in \mathcal{T}, v\in \mathcal{V}_{k,k+t},\\
  &\sum_{j\in \alpha_{k,k+t}^{v}(i)} z_{k,k+t}^{ijv} = \sum_{j\in \beta_{k,k+t}^{v}(i)} z_{k,k+t}^{ijv} &&\forall i\in \mathcal{N},  t\in \mathcal{T}, v\in \mathcal{V}_{k,k+t},\\
  &\sum_{v\in \mathcal{V}_{k,k+t}} \sum_{j\in \alpha_{k,k+t}^{v}(i)} z_{k,k+t}^{ijv} \leq 1 &&\forall i\in \mathcal{N}, t\in \mathcal{T} ,\\
  &u_{k,k+t}^{iv} \leq M\sum_{j\in \alpha_{k,k+t}^{v}(i)} z_{k,k+t}^{ijv} &&\forall i\in \mathcal{N}, t\in \mathcal{T},  v \in \mathcal{V}_{k,k+t},\\
  &\sum_{i\in \mathcal{N}} u_{k,k+t}^{iv} \leq Q &&\forall t\in \mathcal{T}, v\in \mathcal{V}_{k,k+t},\\
  &u_{k,k+t}^{iv}\leq I^{\max} - \hat{I}_{k,k+t}^i &&\forall i\in \mathcal{N}, t\in \mathcal{T},  v \in \mathcal{V}_{k,k+t},\\
  &\hat{I}_{k,k+t+1}^i = \hat{I}_{k,k+t}^i + \sum_{v\in \mathcal{V}_{k,k+t}}u_{k,k+t}^{iv} - \hat{y}_{k,k+t}^i \!\!\!\!\!\!\!\!\!\!\!\!\!\!\!\!\!\!\!&&\forall i\in \mathcal{N},  t\in \mathcal{T} ,\\
  & z_{k,k+t}^{ijv} \in \{0,1\} &&\forall t\in \mathcal{T}, v \in \mathcal{V}_{k,k+t},(i,j)\in E,\\
  &u_{k,k+t}^{iv} \geq 0 &&\forall i\in \mathcal{N}, t\in \mathcal{T},  v\in \mathcal{V}_{k,k+t}.
\end{align}

The penalized objective function is Equation (\ref{mA39}) with Constraints (\ref{mA40}) and (\ref{mA41}).

\begin{align}
\min\ & \sum_{t\in \mathcal{T}}\sum_{(i,j)\in E}\sum_{v\in \mathcal{V}_{k,k+t}}  c_{ij}z_{k,k+t}^{ijv} +\sum_{t\in \mathcal{T}} \sum_{i\in \mathcal{N}}(h (\hat{I}_{k,k+t+1}^i)^+ + e (\hat{I}_{k,k+t+1}^i)^-) + \sum_{i\in \mathcal{N}}\sum_{v\in \mathcal{V}}\lambda_{k}^{iv}(\mu_{k,k}^{iv}+\nu_{k,k}^{iv}), \label{mA39}\\
s.t. \ & \mu_{k,k}^{iv}-\nu_{k,k}^{iv} = u_{k,k+t}^{iv}-\bar{u}_{k,k+t}^{iv} \quad \forall i\in \mathcal{N}, v\in \mathcal{V} ,\label{mA40}\\
& \mu_{k,k}^{iv}, \nu_{k,k}^{iv}\geq 0  \qquad \qquad \qquad \quad \ \forall i\in \mathcal{N}, v\in \mathcal{V}. \label{mA41}
\end{align}

\clearpage
\section{Sensitivity Analysis} \label{6.3}
We show the robustness of ScenPO with regards to varying parameter values for holding and backorder cost parameters. The results are shown in Table \ref{tab:cc0} on a 5 retailer system.
The results show that MQRNN in combination with ScenPO-SS is a robust choice with an average saving of 56.20\% among all instances. 

\begin{table}[hpt]
\centering
\caption{Results for various holding and backorder cost values.}
\label{tab:cc0}
\scriptsize
\begin{tabular}{ c  c  c r r r r r r r r r r r r r}
\hline
\centering
\multirow{2}*{$h$} & \multirow{2}*{$e$} &    \multirow{2}*{Item}    & \multirow{2}*{PI} & \multirow{2}*{EV} & \multirow{2}*{EMP} & \multicolumn{3}{c}{Prediction-Focused Learning}         & \multicolumn{3}{c}{ScenPO-SS} &  \multicolumn{3}{c}{ScenPO-TS} \\
\cmidrule(r){7-9}
\cmidrule(r){10-12}
\cmidrule(r){13-15}
     &    &              &            &           &              & MLE         & LSTM    & MQRNN   & MLE              & LSTM    & MQRNN   & MLE              & LSTM     & MQRNN   \\
\hline
\multirow{9}*{0.15}  & \multirow{9}*{1.5}  &Transportation cost/Cost & 75\%  & 61\% & 64\%  & 61\%  & 66\% & 66\% & 58\% & 63\% & 62\% & 73\%  & 70\% & 72\% \\
     &    & Holding cost/Cost      & 25\%  & 25\% & 29\%  & 26\%  & 25\% & 23\% & 38\% & 28\% & 31\% & 18\%  & 14\% & 15\% \\
     &    & Backorder cost/Cost     & 0\%   & 14\% & 7\%   & 13\%  & 9\%  & 11\% & 4\%  & 9\%  & 6\%  & 9\%   & 16\% & 13\% \\
   &      & Cost                & 684   & 866  & 890   & 920   & 822  & 831  & 802  & 783  & 775  & 891   & 880  & 863  \\
      &   & $\Delta$Cost               & 0    & 182  & 206   & 236   & 138  & 147  & 118  & 99   & 91   & 207   & 196  & 179  \\
    &     & Saving              & 100\%   & 0\%    & -13\% & -29\% & 24\% & 19\% & 35\% & 45\% & 50\% & -14\% & -8\% & 2\%  \\
\hline
\end{tabular}
\end{table}

\clearpage

\section{Complete Experimental Results} \label{a2}

\begin{table}[hpt]
\centering
\caption{Detailed numerical results of ScenPO and the proposed benchmarks for a 5-retailer system.}
\label{tab:gd5}
\scriptsize
\begin{tabular}{ c  c  c r r r r r r r r r r r r}
\hline
\centering
  &  & & \multirow{2}*{PI} & \multirow{2}*{EV} & \multirow{2}*{EMP} & \multicolumn{3}{c}{Prediction-Focused Learning}         & \multicolumn{3}{c}{ScenPO-SS} &  \multicolumn{3}{c}{ScenPO-TS} \\
\cmidrule(r){7-9}
\cmidrule(r){10-12}
\cmidrule(r){13-15}
 &      &              &            &           &              & MLE         & LSTM    & MQRNN   & MLE              & LSTM    & MQRNN   & MLE              & LSTM     & MQRNN   \\
\hline
 \multirow{18}*{5} & \multirow{6}*{\textit{Random}} & Cost    & 718   & 1183 & 1207  & 1284  & 1166 & 1152 & 1013 & 998  & 992  & 1149 & 1169 & 1150 \\
   &        & $\Delta$Cost   & 0     & 465  & 489   & 566   & 448  & 434  & 295  & 280  & 274  & 431  & 451  & 432  \\
   &        & Saving  &100\%     & 0\%  & -5\%  & -22\% & 4\%  & 6\%  & 36\% & 40\% & 41\% & 7\%  & 3\%  & 7\%  \\
   &        & MSE     & 0     & 289  & 638   & 628   & 302  & 299  &      &      &      &      &      &      \\
   &        & Service & 100\% & 79\% & 84\%  & 77\%  & 79\% & 80\% & 90\% & 89\% & 89\% & 80\% & 77\% & 79\% \\
   &        & Time (s)  & 5     & 6    & 5     & 1     & 6    & 6    & 9    & 8    & 8    & 1109 & 1095 & 768  \\
\hline
\end{tabular}
\end{table}

\begin{table}[hpt]
\centering
\caption{Detailed numerical results of ScenPO and the proposed benchmarks for a 7-retailer system.}
\label{tab:gd7}
\scriptsize
\begin{tabular}{ c  c c  r r r r r r r r r r r r}
\hline
\centering
 &  &   & \multirow{2}*{PI} & \multirow{2}*{EV} & \multirow{2}*{EMP} & \multicolumn{3}{c}{Prediction-Focused Learning}         & \multicolumn{3}{c}{ScenPO-SS} &  \multicolumn{3}{c}{ScenPO-TS} \\
\cmidrule(r){7-9}
\cmidrule(r){10-12}
\cmidrule(r){13-15}
&      &              &            &           &              & MLE         & LSTM    & MQRNN   & MLE              & LSTM    & MQRNN   & MLE              & LSTM     & MQRNN   \\
\hline
\multirow{18}*{5} & \multirow{6}*{\textit{Random}} & Cost    & 977   & 1634 & 1746  & 1749  & 1639 & 1654 & 1362 & 1400 & 1368 & 1576 & 1597 & 1570 \\
   &        & $\Delta$Cost   & 0  & 657  & 769   & 772   & 662  & 677  & 385  & 424  & 391  & 599  & 620  & 593  \\
   &        & Saving  &100\%    & 0\%  & -17\% & -18\% & -1\% & -3\% & 41\% & 36\% & 40\% & 9\%  & 6\%  & 10\% \\
   &        & MSE     & 0     & 292  & 612   & 627   & 309  & 306  &      &      &      &      &      &      \\
   &        & Service & 100\% & 77\% & 79\%  & 77\%  & 78\% & 76\% & 91\% & 89\% & 89\% & 78\% & 76\% & 78\% \\
   &        & Time (s)  & 14    & 16   & 14    & 3     & 15   & 15   & 35   & 31   & 30   & 1479 & 1466 & 1298 \\
\hline
\end{tabular}
\end{table}

\begin{table}[hpt]
\centering
\caption{Detailed numerical results of ScenPO and the proposed benchmarks for a 10-retailer system.}
\label{tab:gd10}
\scriptsize
\begin{tabular}{c c  c r r r r r r r r r r r r}
\hline
\centering
&  &  & \multirow{2}*{PI} & \multirow{2}*{EV} & \multirow{2}*{EMP} & \multicolumn{3}{c}{Prediction-Focused Learning}         & \multicolumn{3}{c}{ScenPO-SS} &  \multicolumn{3}{c}{ScenPO-TS} \\
\cmidrule(r){7-9}
\cmidrule(r){10-12}
\cmidrule(r){13-15}
&      &              &            &           &              & MLE         & LSTM    & MQRNN   & MLE              & LSTM    & MQRNN   & MLE              & LSTM     & MQRNN   \\
\hline
\multirow{18}*{5} & \multirow{6}*{\textit{Random}} & Cost    & 1290  & 2302 & 2454  & 2662  & 2325 & 2358 & 1844 & 1877 & 1849 & 2152 & 2226 & 2143 \\
   &        & $\Delta$Cost   & 0   & 1012 & 1164  & 1372  & 1035 & 1068 & 553  & 587  & 559  & 862  & 936  & 853  \\
   &        & Saving  & 100\%    & 0\%  & -15\% & -36\% & -2\% & -6\% & 45\% & 42\% & 45\% & 15\% & 8\%  & 16\% \\
   &        & MSE     & 0     & 277  & 570   & 613   & 294  & 290  &      &      &      &      &      &      \\
   &        & Service & 100\% & 70\% & 75\%  & 69\%  & 70\% & 68\% & 91\% & 89\% & 88\% & 73\% & 72\% & 76\% \\
   &        & Time (s)  & 103   & 110  & 88    & 16    & 97   & 115  & 690  & 542  & 615  & 4587 & 3977 & 3833 \\
\hline
\end{tabular}
\end{table}

\begin{table}[hpt]
\centering
\caption{Comparison of the various scenario generation approaches for a 5-retailer system.}
\scriptsize
\label{tab:expend5}
\begin{tabular}{cccrrrrrrrrr}
\hline
 & &   & \multirow{2}*{PI} & \multirow{2}*{EV} & \multicolumn{5}{c}{residuals-based}        & \multicolumn{2}{c}{quantile-based}            \\
\cmidrule(r){6-10}
\cmidrule(r){11-12}
         &           &         &            &           &  OLS-ER & OLS-J & OLS-JP      & ES & NN & QuantReg          & ScenPO \\
\hline
\multirow{15}{*}{5}  & \multirow{5}{*}{\textit{Random}} & Cost         & 718   & 1183 & 998             & 989   & 989    & 1093 & 1026 & 1005            & 992    \\
                     &                                                   & $\Delta$Cost & 0     & 465  & 280             & 271   & 271    & 375  & 308  & 287             & 274    \\
                     &                                                   & Saving       & 100\% & 0\%  & 40\%            & 42\%  & 42\%   & 19\% & 34\% & 38\%            & 41\%   \\
                     &                                                   & Service      & 100\% & 79\% & 93\%            & 93\%  & 93\%   & 91\% & 90\% & 89\%            & 89\%   \\
                     &                                                   & Time (s)     & 5     & 6    & 8               & 8     & 9      & 7    & 8    & 8               & 8      \\
                     \cline{2-12}
                     & \multirow{5}{*}{\textit{Trend}}  & Cost         & 691   & 1047 & 778             & 769   & 768    & 762  & 849  & 775             & 774    \\
\hline
\end{tabular}
\end{table}

\begin{table}[hpt]
\centering
\caption{Comparison of the various scenario generation approaches for a 7-retailer system.}
\scriptsize
\label{tab:expend7}
\begin{tabular}{cccrrrrrrrrr}
\hline
 & &   & \multirow{2}*{PI} & \multirow{2}*{EV} & \multicolumn{5}{c}{residuals-based}        & \multicolumn{2}{c}{quantile-based}            \\
\cmidrule(r){6-10}
\cmidrule(r){11-12}
         &           &         &            &           &  OLS-ER & OLS-J & OLS-JP  &ES    & NN & QuantReg          & ScenPO \\
\hline
\multirow{15}{*}{5}  & \multirow{5}{*}{\textit{Random}} & Cost         & 977   & 1634 & 1384            & 1380  & 1382   & 1509 & 1381 & 1373            & 1368   \\
                     &                                                   & $\Delta$Cost & 0     & 657  & 407             & 403   & 405    & 532  & 404  & 396             & 391    \\
                     &                                                   & Saving       & 100\% & 0\%  & 38\%            & 39\%  & 38\%   & 19\% & 38\% & 40\%            & 40\%   \\
                     &                                                   & Service      & 100\% & 77\% & 91\%            & 92\%  & 92\%   & 91\% & 90\% & 89\%            & 89\%   \\
                     &                                                   & Time (s)     & 14    & 16   & 36              & 39    & 51     & 33   & 38   & 31              & 30     \\
                     \cline{2-12}
                     & \multirow{5}{*}{\textit{Trend}}  & Cost         & 902   & 1418 & 1009            & 1012  & 1010   & 987  & 1111 & 1024            & 1022   \\
\hline
\end{tabular}
\end{table}

\begin{table}[hpt]
\centering
\caption{Comparison of the various scenario generation approaches for a 10-retailer system.}
\scriptsize
\label{tab:expend10}
\begin{tabular}{cccrrrrrrrrr}
\hline
 & &   & \multirow{2}*{PI} & \multirow{2}*{EV} & \multicolumn{5}{c}{residuals-based}        & \multicolumn{2}{c}{quantile-based}            \\
\cmidrule(r){6-10}
\cmidrule(r){11-12}
         &           &         &            &           &  OLS-ER & OLS-J & OLS-JP      & ES & NN & QuantReg          & ScenPO \\
\hline
\multirow{15}{*}{5}  & \multirow{5}{*}{\textit{Random}} & Cost         & 1290  & 2302 & 1858            & 1857  & 1850   & 2064 & 1879 & 1850            & 1849   \\
                     &                                                   & $\Delta$Cost & 0     & 1012 & 568             & 566   & 559    & 773  & 588  & 560             & 559    \\
                     &                                                   & Saving       & 100\% & 0\%  & 44\%            & 44\%  & 45\%   & 24\% & 42\% & 45\%            & 45\%   \\
                     &                                                   & Service      & 100\% & 70\% & 91\%            & 93\%  & 94\%   & 92\% & 89\% & 88\%            & 88\%   \\
                     &                                                   & Time (s)     & 103   & 110  & 752             & 690   & 671    & 745  & 606  & 659             & 615    \\
                     \cline{2-12}
                     & \multirow{5}{*}{\textit{Trend}}  & Cost         & 1182  & 2001 & 1335            & 1344  & 1339   & 1317 & 1462 & 1340            & 1341   \\
\hline
\end{tabular}
\end{table}

\begin{table}[hpt]
\centering
\caption{Detailed numerical results of CombOptNet for a 5-retailer system.}
\label{tab:comb}
\scriptsize
\begin{tabular}{c c  c r r r r r r r r r r r r r}
\hline
\centering
&   & \multirow{2}*{PI} & \multirow{2}*{EV} & \multirow{2}*{EMP} & \multicolumn{3}{c}{Prediction-Focused Learning}         & \multicolumn{3}{c}{ScenPO-SS} &  \multicolumn{3}{c}{ScenPO-TS}&  \multirow{2}*{CombOptNet}\\
\cmidrule(r){6-8}
\cmidrule(r){9-11}
\cmidrule(r){12-14}
                                 &              &       &      &       & MLE                        & LSTM & MQRNN & MLE       & LSTM   & MQRNN & MLE       & LSTM  & MQRNN &            \\
                                 \hline
\multirow{5}*{\textit{Random}} & Cost         & 718   & 1183 & 1207  & 1284                       & 1166 & 1152  & 1013      & 998    & 992   & 1149      & 1169  & 1150  & 1234       \\
                                 & $\Delta$Cost & 0     & 465  & 489   & 566                        & 448  & 434   & 295       & 280    & 274   & 431       & 451   & 432   & 516        \\
                                 & Saving       & 100\% & 0\%  & -5\%  & -22\%                      & 4\%  & 6\%   & 36\%      & 40\%   & 41\%  & 7\%       & 3\%   & 7\%   & -11\%      \\
                                 & Service      & 100\% & 79\% & 84\%  & 77\%                       & 79\% & 80\%  & 90\%      & 89\%   & 89\%  & 80\%      & 77\%  & 79\%  & 96\%       \\
                                 & Time (s)     & 5     & 6    & 5     & 1                          & 6    & 6     & 9         & 8      & 8     & 1109      & 1095  & 768   & 5          \\
                                 \hline
\multirow{5}*{\textit{Trend}}  & Cost         & 691   & 1047 & 1042  & 1099                       & 886  & 873   & 938       & 915    & 774   & 966       & 987   & 906   & 1051       \\
\hline
\end{tabular}
\end{table}

\begin{table}
\centering
\caption{Detailed numerical results of ScenPO and the proposed benchmarks for real-life Instance 1.}
\label{tab:rd1}
\scriptsize
\begin{tabular}{c c c r r r r r r r r r r r r}
\hline
\centering
\multirow{2}*{\#Retailer} & \multirow{2}*{$\ell$} & \multirow{2}*{Item}         & \multirow{2}*{PI} & \multirow{2}*{EV} & \multirow{2}*{EMP} & \multicolumn{3}{c}{Prediction-Focused Learning}         & \multicolumn{3}{c}{ScenPO-SS} &  \multicolumn{3}{c}{ScenPO-TS} \\
\cmidrule(r){7-9}
\cmidrule(r){10-12}
\cmidrule(r){13-15}
           &          &              &            &           &              & MLE         & LSTM    & MQRNN   & MLE              & LSTM    & MQRNN   & MLE              & LSTM     & MQRNN   \\
\hline
\multirow{18}*{5}  & \multirow{6}*{5}  & Cost    & 413   & 2198  & 2390  & 2513  & 1519  & 1788  & 3011  & 4445   & 1464  & 2278  & 3407  & 2241  \\
   &    & $\Delta$Cost   & 0     & 1784  & 1977  & 2100  & 1106  & 1375  & 2598  & 4032   & 1051  & 1865  & 2993  & 1827  \\
   &    & Saving  & 100\%     & 0\%   & -11\% & -18\% & 38\%  & 23\%  & -46\% & -126\% & 41\%  & -5\%  & -68\% & -2\%  \\
   &    & MSE     & 0     & 58    & 112   & 282   & 33    & 39    &       &        &       &       &       &       \\
   &    & Service & 97\%  & 77\%  & 74\%  & 92\%  & 85\%  & 75\%  & 99\%  & 100\%  & 95\%  & 91\%  & 94\%  & 74\%  \\
   &    & Time (s)  & 2     & 0     & 0     & 0     & 1     & 1     & 2     & 1      & 3     & 792   & 778   & 653   \\
\hline
\end{tabular}
\end{table}

\begin{table}
\centering
\caption{Detailed numerical results of ScenPO and the proposed benchmarks for real-life Instance 2.}
\label{tab:rd2}
\scriptsize
\begin{tabular}{c c c r r r r r r r r r r r r}
\hline
\centering
\multirow{2}*{\#Retailer} & \multirow{2}*{$\ell$} & \multirow{2}*{Item}         & \multirow{2}*{PI} & \multirow{2}*{EV} & \multirow{2}*{EMP} & \multicolumn{3}{c}{Prediction-Focused Learning}         & \multicolumn{3}{c}{ScenPO-SS} &  \multicolumn{3}{c}{ScenPO-TS} \\
\cmidrule(r){7-9}
\cmidrule(r){10-12}
\cmidrule(r){13-15}
           &          &              &            &           &              & MLE         & LSTM    & MQRNN   & MLE              & LSTM    & MQRNN   & MLE              & LSTM     & MQRNN   \\
\hline
\multirow{18}*{5}  & \multirow{6}*{5}  & Cost    & 3980  & 8554  & 9320  & 10345 & 7257  & 8209  & 7514  & 7955   & 6436  & 9897  & 11410 & 9521  \\
   &    & $\Delta$Cost   & 0    & 4574  & 5340  & 6365  & 3277  & 4229  & 3534  & 3976   & 2456  & 5917  & 7430  & 5542  \\
   &    & Saving  & 100\%   & 0\%   & -17\% & -39\% & 28\%  & 8\%   & 23\%  & 13\%   & 46\%  & -29\% & -62\% & -21\% \\
   &    & MSE     & 0     & 300   & 514   & 465   & 179   & 175   &       &        &       &       &       &       \\
   &    & Service & 98\%  & 71\%  & 71\%  & 69\%  & 79\%  & 67\%  & 85\%  & 96\%   & 88\%  & 67\%  & 65\%  & 60\%  \\
   &    & Time (s)  & 7     & 4     & 1     & 0     & 6     & 6     & 5     & 2      & 5     & 1262  & 927   & 1066  \\
\hline
\end{tabular}
\end{table}

\begin{table}[hpt]
\centering
\caption{Comparison of the various scenario generation approaches for real-life Instance 1.}
\label{tab:expandrd1}
\scriptsize
\begin{tabular}{cccrrrrrrrrr}
\hline
\multirow{2}*{\#Retailer} & \multirow{2}*{$\ell$} & \multirow{2}*{Item}   & \multirow{2}*{PI} & \multirow{2}*{EV} & \multicolumn{5}{c}{residuals-based}        & \multicolumn{2}{c}{quantile-based}            \\
\cmidrule(r){6-10}
\cmidrule(r){11-12}
         &           &         &            &           & OLS-ER & OLS-J & OLS-JP   & ES   & NN & QuantReg          & ScenPO \\
\hline
\multirow{15}*{5}          &\multirow{5}*{ 5}  & Cost         & 413   & 2198    & 3292            & 2211  & 2208   & 3002  & 2224  & 1615            & 1464   \\
                     &                     & $\Delta$Cost & 0     & 1784    & 2879            & 1798  & 1795   & 2589  & 1810  & 1201            & 1051   \\
                     &                     & Saving       & 100\% & 0\%     & -61\%           & -1\%  & -1\%   & -45\% & -1\%  & 33\%            & 41\%   \\
                     &                     & Service      & 97\%  & 77\%    & 99\%            & 96\%  & 96\%   & 99\%  & 99\%  & 93\%            & 95\%   \\
                     &                     & Time (s)     & 2     & 0       & 2               & 2     & 2      & 2     & 3     & 2               & 3      \\
                     \cline{2-12}
                     & \multirow{5}{*}{7}  & Cost         & 349   & 2097    & 3034            & 2432  & 2434   & 2872  & 2305  & 1472            & 1377   \\
\hline
\end{tabular}
\end{table}

\begin{table}[hpt]
\centering
\caption{Comparison of the various scenario generation approaches for real-life Instance 2.}
\label{tab:expandrd2}
\scriptsize
\begin{tabular}{cccrrrrrrrrr}
\hline
\multirow{2}*{\#Retailer} & \multirow{2}*{$\ell$} & \multirow{2}*{Item}   & \multirow{2}*{PI} & \multirow{2}*{EV} & \multicolumn{5}{c}{residuals-based}        & \multicolumn{2}{c}{quantile-based}            \\
\cmidrule(r){6-10}
\cmidrule(r){11-12}
         &           &         &            &           & OLS-ER & OLS-J & OLS-JP   &ES   & NN & QuantReg          & ScenPO \\
\hline
\multirow{15}{*}{5}  & \multirow{5}{*}{5}  & Cost         & 3980  & 8554  & 7050            & 6827  & 6588   & 6914  & 6625  & 6689            & 6436   \\
                     &                     & $\Delta$Cost & 0     & 4574  & 3070            & 2847  & 2608   & 2934  & 2646  & 2709            & 2456   \\
                     &                     & Saving       & 100\% & 0\%   & 33\%            & 38\%  & 43\%   & 36\%  & 42\%  & 41\%            & 46\%   \\
                     &                     & Service      & 98\%  & 71\%  & 89\%            & 85\%  & 84\%   & 89\%  & 86\%  & 83\%            & 88\%   \\
                     &                     & Time (s)     & 7     & 4     & 3               & 5     & 6      & 5     & 6     & 2               & 5      \\
                     \cline{2-12}
                     & \multirow{5}{*}{7}  & Cost         & 3970  & 8416  & 6742            & 7027  & 7346   & 6991  & 6392  & 6643            & 6125   \\
\hline
\end{tabular}
\end{table}

\end{APPENDICES}		
\end{document}